\documentclass[12pt]{article}
\usepackage{amssymb}
\usepackage{amsmath}
\usepackage{amsthm}
\usepackage{enumitem}
\usepackage{graphicx}
\begin{document}

\newtheorem{theorem}{Theorem}
\newtheorem{corollary}[theorem]{Corollary}
\newtheorem{lemma}[theorem]{Lemma}
\theoremstyle{definition}
\newtheorem{definition}{Definition}

\title{Mixed-Mode Oscillations in a Stochastic, Piecewise-Linear System}
\author{
D.J.W.~Simpson and R.~Kuske\thanks{
The authors acknowledge support from an NSERC Discovery Grant.
}\\
Department of Mathematics\\
University of British Columbia\\
Vancouver, BC, V6T1Z2}
\maketitle

\begin{abstract}
We analyze a piecewise-linear FitzHugh-Nagumo model.
The system exhibits a canard near which
both small amplitude and large amplitude periodic orbits exist.
The addition of small noise induces mixed-mode oscillations (MMOs)
in the vicinity of the canard point.
We determine the effect of each model parameter on the
stochastically driven MMOs.
In particular we show that
any parameter variation (such as a modification of the piecewise-linear
function in the model) that leaves the ratio of noise
amplitude to time-scale separation unchanged typically has little effect
on the width of the interval of the primary bifurcation parameter
over which MMOs occur.
In that sense, the MMOs are robust.
Furthermore we show that the piecewise-linear model
exhibits MMOs more readily than the classical FitzHugh-Nagumo
model for which a cubic polynomial is the only nonlinearity.
By studying a piecewise-linear model we are able to
explain results using analytical expressions
and compare these with numerical investigations.
\end{abstract}

\section{Introduction}
\label{sec:INTRO}

\begin{figure}[b!]
\begin{center}
\includegraphics[height=4cm]{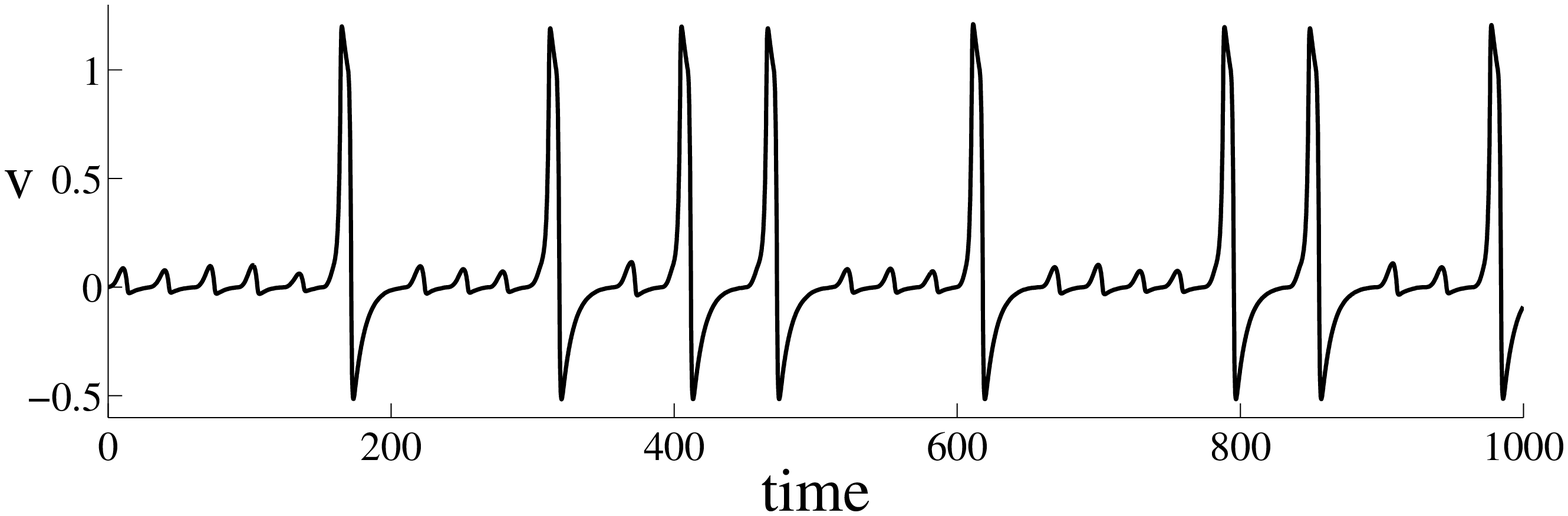}
\caption{
A time series illustrating MMOs exhibited by (\ref{eq:fhn}) with (\ref{eq:fpwl}).
The parameter values are the same as in Fig.~\ref{fig:ssEx}-A.
\label{fig:timeSeries}
}
\end{center}
\end{figure}

Oscillatory dynamics involving oscillations with greatly differing
amplitudes, known as mixed-mode oscillations (MMOs), see Fig.~\ref{fig:timeSeries},
are important in neuron models \cite{ErMc08}
and in a multitude of chemical reactions \cite{Ba88,PeSc92},
refer to \cite{DeGu10} for a recent review.
Yet there are many open questions regarding the
creation, robustness and bifurcations of MMOs.
A variety of mechanisms generate MMOs in deterministic systems.
Alternatively MMOs may be noise-induced;
there are also several scenarios by which this may occur.

We study the following form of the FitzHugh-Nagumo (FHN) model
with small, additive, white noise:
\begin{equation}
\begin{split}
dv & = (f(v) - w)~dt \;, \\
dw & = \varepsilon ( \alpha v - \sigma w - \lambda )~dt + D~dW \;,
\end{split}
\label{eq:fhn}
\end{equation}
where $v$ represents a potential,
$w$ is a recovery variable
and $W$ is a standard Brownian motion.
The FHN model is used as a prototypical model of excitable dynamics
in a range of scientific fields \cite{RoGe00,KeSn98}.
Here $\alpha$ is a positive constant
and $\lambda \in \mathbb{R}$,
which is regarded as the main bifurcation parameter,
controls the growth of oscillations, as seen below.
The small parameter $\varepsilon \ll 1$,
represents the time-scale separation and
$D \ll 1$ is the noise amplitude ($\varepsilon,D > 0$).
Values of $\varepsilon$ and $D$ used in, for instance \cite{MuVa08,MaNe01},
are no larger than the values considered here.
By scaling we may assume $\sigma = 1$,
except in the special case $\sigma = 0$
which corresponds to the van der Pol model
(and in this case we may further assume $\alpha = 1$).
We assume that $f : \mathbb{R} \to \mathbb{R}$
is continuous and roughly of cubic shape.
For simplicity, we assume that $f$ has a local minimum at $(0,0)$
and a local maximum at $(1,1)$,
regardless of the precise function chosen.

\begin{figure}[t!]
\begin{center}
\setlength{\unitlength}{1cm}
\begin{picture}(14.7,6.4)
\put(0,0){\includegraphics[width=7.2cm,height=6cm]{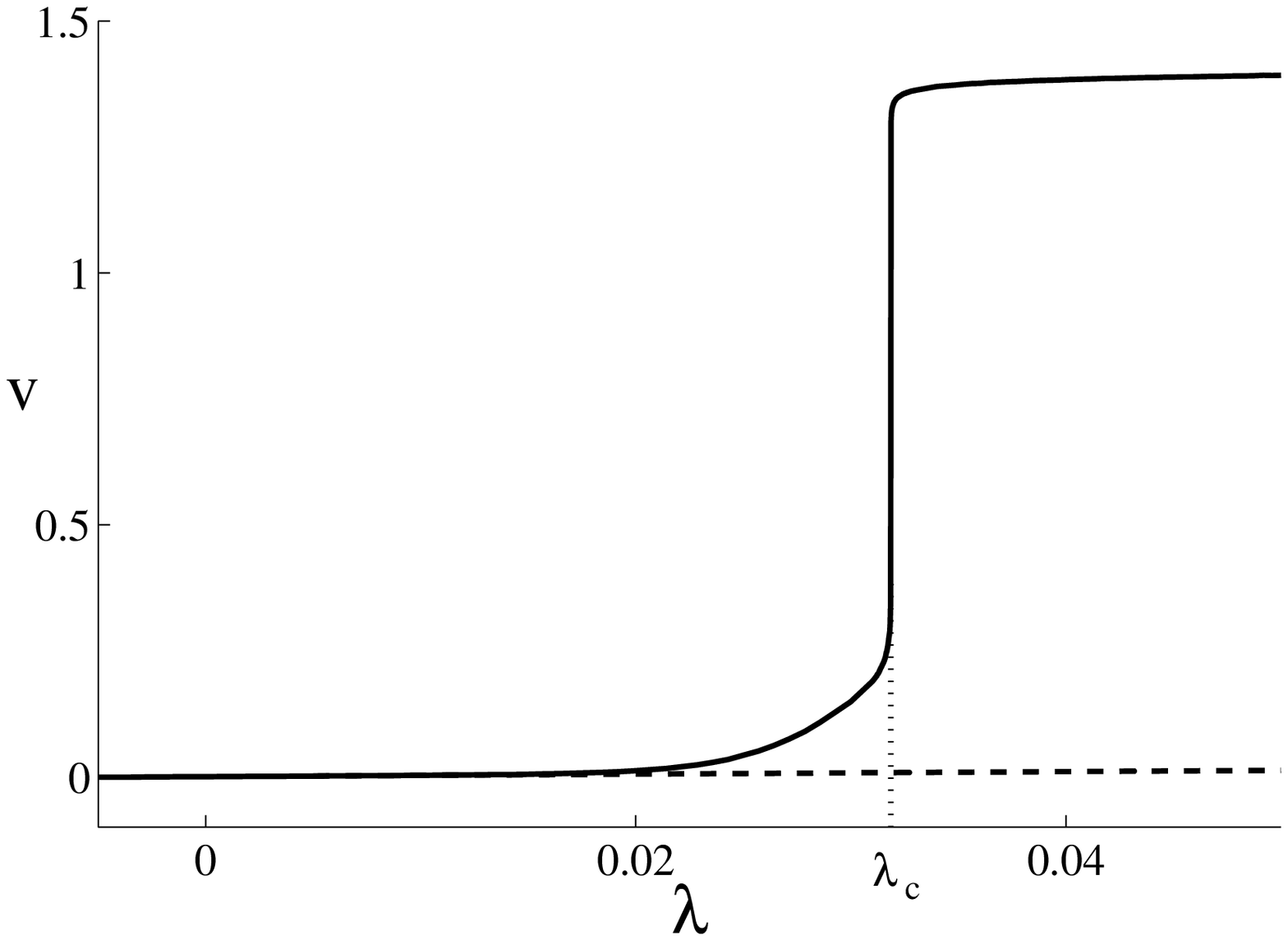}}
\put(7.7,0){\includegraphics[width=7.2cm,height=6cm]{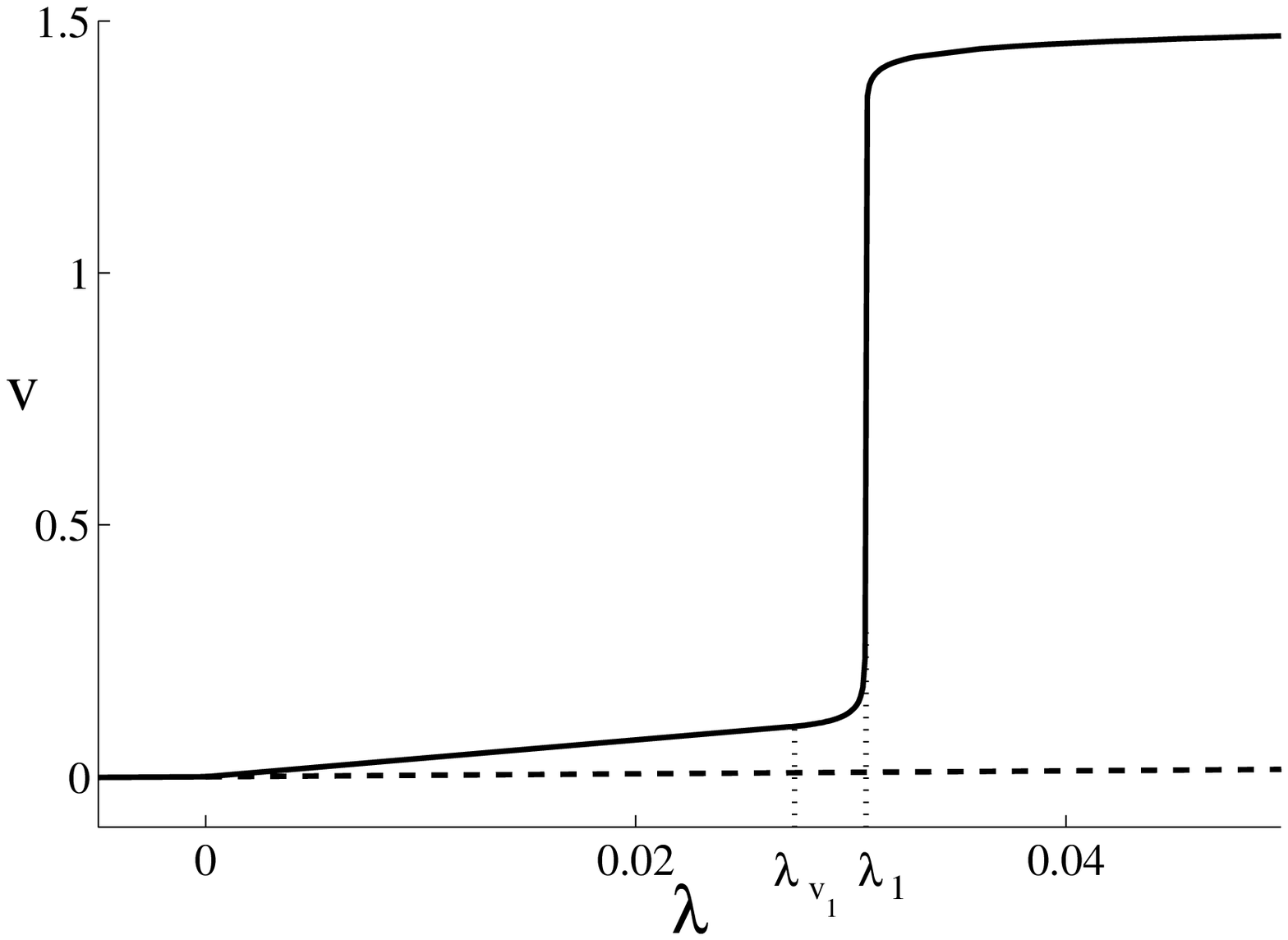}}
\put(1.2,6){\large \sf \bfseries A}
\put(8.9,6){\large \sf \bfseries B}
\put(3.4,6.3){smooth}
\put(11.1,6.3){PWL}
\end{picture}
\caption{
Bifurcation diagrams of (\ref{eq:fhn}) in the absence of noise
(i.e.~$D = 0$) with (\ref{eq:fcubic})
in panel A and with (\ref{eq:fpwl}) in panel B.
In each panel the solid curve for $\lambda > 0$
corresponds to the maximum $v$-value of a stable periodic orbit;
the remaining curves correspond to the equilibrium
which is unstable for $\lambda > 0$.
In panel A a canard explosion occurs near the canard point, $\lambda_c$;
in panel B a canard explosion occurs near $\lambda_1$
at which point the stable periodic orbit has a maximum value of $1$.
The parameter values used are
$\varepsilon = 0.04$,
$(\alpha,\sigma) = (4,1)$,
$(\eta_L,\eta_R) = (-2,-1)$ and
$(v_1,w_1) = (0.1,0.05)$.
\label{fig:canardBifDiag}
}
\end{center}
\end{figure}

\begin{figure}[ht]
\begin{center}
\setlength{\unitlength}{1cm}
\begin{picture}(14.7,6.4)
\put(0,0){\includegraphics[width=7.2cm,height=6cm]{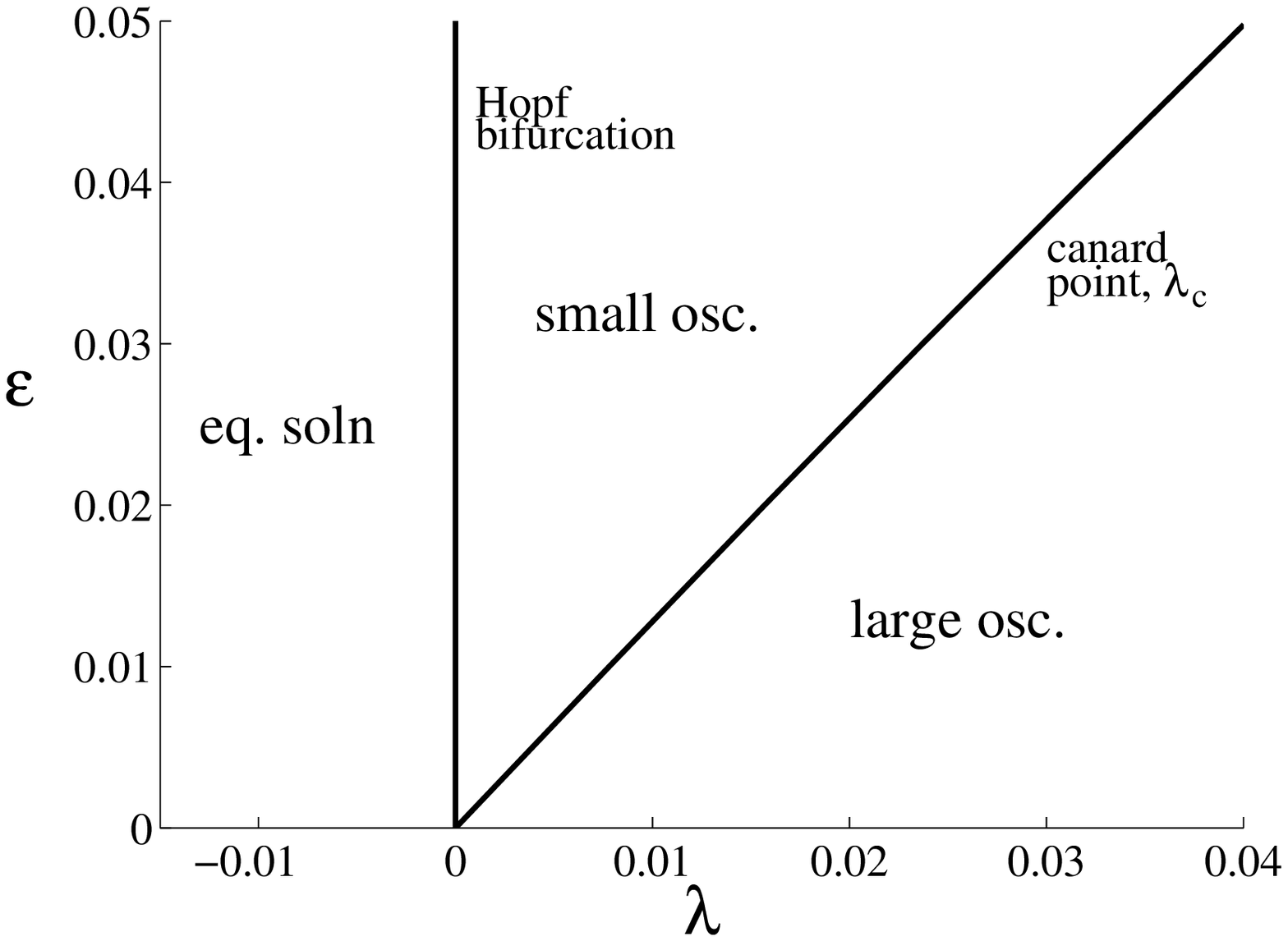}}
\put(7.7,0){\includegraphics[width=7.2cm,height=6cm]{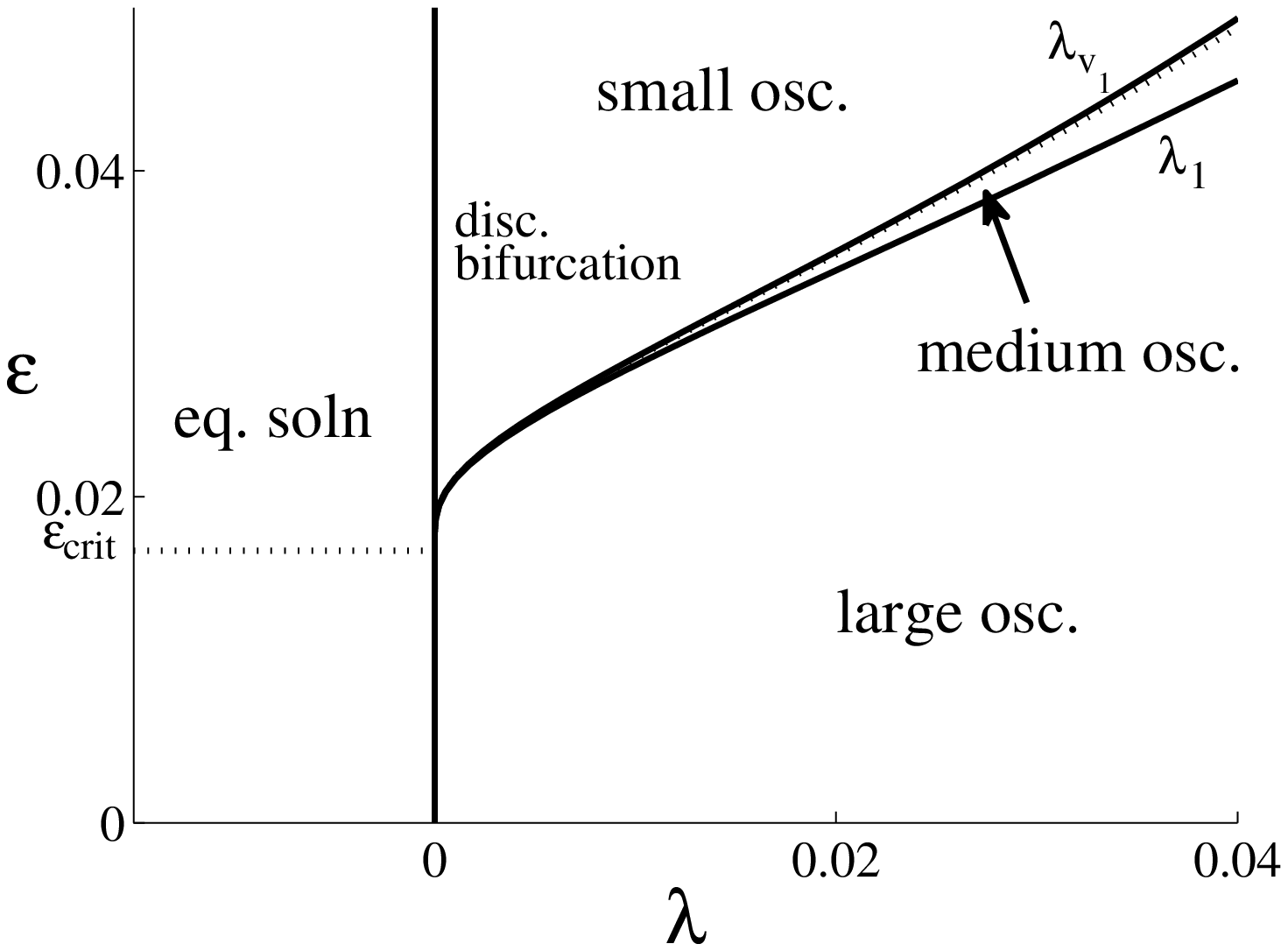}}
\put(1.2,6){\large \sf \bfseries A}
\put(8.9,6){\large \sf \bfseries B}
\put(3.4,6.3){smooth}
\put(11.1,6.3){PWL}
\end{picture}
\caption{
Two parameter bifurcation diagrams of the smooth and PWL versions
of (\ref{eq:fhn}) with the same parameter values as in
Fig.~\ref{fig:canardBifDiag}.
The smooth system has a well-defined canard point, $\lambda_c$
\cite{KrSz01,KrSz01b}, whereas for the PWL system
we consider the two values, $\lambda_{v_1}$ and $\lambda_1$,
described in the text.
In both panels we have indicated the attracting solution for each
region bounded by the solid curves.
The dotted curve in panel B corresponds to the approximation
(\ref{eq:lamApprox}) derived below;
$\varepsilon_{\rm crit}$ is given by (\ref{eq:eCrit}).
Note that in contrast to the remainder of this paper,
in panel A the distinction between small and large oscillations
is determined by $\lambda_c$ and not (\ref{eq:size}).
\label{fig:lamEe}
}
\end{center}
\end{figure}

If $f$ is a cubic,
as originally taken by FitzHugh \cite{Fi61}
and Nagumo {\em et.~al.}~\cite{NaAr62},
then, by the above requirements, the cubic must be
\begin{equation}
f(v) = 3v^2 - 2v^3 \;.
\label{eq:fcubic}
\end{equation}
Fig.~\ref{fig:canardBifDiag}-A illustrates the
role of the parameter $\lambda$
for (\ref{eq:fhn}) with (\ref{eq:fcubic}) in the absence of noise.
A small amplitude periodic orbit is created in
a Hopf bifurcation at $\lambda = 0$.
For the parameters used in Fig.~\ref{fig:canardBifDiag},
this periodic orbit is stable and its amplitude increases with $\lambda$.
Near $\lambda_c$ the amplitude increases exponentially.
This rapid growth is known as a canard explosion
and is due to time-scale separation and global dynamics
\cite{BeCa81,Ec83,BaEr86,BaEr92}.
The value of the {\em canard point}, $\lambda_c$,
which is well-defined for smooth systems \cite{KrSz01,KrSz01b},
decreases to zero with $\varepsilon$,
as shown in Fig.~\ref{fig:lamEe}-A.
Over an order $\varepsilon$ range of $\lambda$ values,
(\ref{eq:fhn}) with (\ref{eq:fcubic})
may either settle to equilibrium, exhibit small amplitude oscillations,
or exhibit large amplitude oscillations (relaxation oscillations).

As in \cite{Co08,ToGe03}, here we study a piecewise-linear (PWL)
FHN model so that, in the presence of noise, the system is amenable to
a rigorous analysis without the need for an approximation
or limiting scenario.
PWL models are commonly used in circuit systems
\cite{BaVe01,ZhMo03,Ts03}.
A PWL version of a driven van der Pol oscillator is studied in \cite{SeIn04}
to explain the breakdown of canards in experiments.
We consider the continuous, PWL function
\begin{equation}
f(v) = \left\{ \begin{array}{lc}
\eta_L v \;,& v \le 0 \\
\eta_1 v \;,& 0 < v \le v_1 \\
\eta_2 (v-v_1) + w_1 \;,& v_1 < v \le 1 \\
\eta_R (v-1) + 1 \;,& v > 1
\end{array} \right. \;,
\label{eq:fpwl}
\end{equation}
where $0 < v_1, w_1 < 1$, $\eta_L, \eta_R < 0$, and
\begin{equation}
\eta_1 = \frac{w_1}{v_1} \;, \qquad
\eta_2 = \frac{1-w_1}{1-v_1} \;.
\label{eq:eta1eta2}
\end{equation}
We state the particular form here in order to briefly illustrate key differences
between the smooth and PWL FHN models.
Further motivation for (\ref{eq:fpwl}) is given in \S\ref{sec:DETER}
and shown in Fig.~\ref{fig:pp1},
As shown in Fig.~\ref{fig:canardBifDiag}-B,
(\ref{eq:fhn}) with (\ref{eq:fpwl}) may exhibit a canard explosion.
The canard point, $\lambda_c$, is not well-defined for
this system because it lacks global differentiability.
Instead we consider the values, $\lambda_{v_1}$ and $\lambda_1$,
at which the maximum $v$-value of the periodic orbit of 
(\ref{eq:fhn}) with (\ref{eq:fpwl}) in the absence of noise
is $v_1$ and $1$ respectively.
The piecewise nature of (\ref{eq:fpwl}) leads to a natural
classification of periodic orbits and oscillations
of (\ref{eq:fhn}) with (\ref{eq:fpwl}).
(Typically we refer to one complete revolution about the equilbrium
as a single oscillation.)
With Fig.~\ref{fig:canardBifDiag}-B in mind,
if $v_{\rm max}$ is the maximum $v$-value
of a periodic orbit or single oscillation
we declare that the orbit or oscillation is
\begin{equation}
\begin{split}
{\rm small~if~} & 0 < v_{\rm max} \le v_1 \;, \\
{\rm medium~if~} & v_1 < v_{\rm max} \le 1 \;, \\
{\rm large~if~} & v_{\rm max} > 1 \;.
\end{split}
\label{eq:size}
\end{equation}
Fig.~\ref{fig:lamEe}-B illustrates typical dependence of 
$\lambda_{v_1}$ and $\lambda_1$ on $\varepsilon$.
In particular we notice that for
a fixed choice of the slopes, $\eta_j$, in (\ref{eq:fpwl}),
the PWL version of the FHN model does not exhibit
small oscillations for arbitrarily small $\varepsilon$.
This is because the two eigenvalues associated with the
equilibrium for small $\lambda > 0$ are real-valued
for sufficiently small $\varepsilon$ negating the possibility of
small oscillations, see \S\ref{sec:DETER}.
Unlike for the van der Pol model,
values of $\varepsilon$ that are relevant for the
FHN model are usually sufficiently large for small oscillations to
be important in the PWL model.

The effect of noise in (\ref{eq:fhn})
has seen significant recent attention, see for instance \cite{LiGa04,LeVa05,PiKu97}.
Noise may induce regular oscillations in (\ref{eq:fhn})
when in the absence of noise there are no oscillations.
There is more than one mechanism that may cause this, most notably
{\em stochastic resonance} \cite{GaHa98}
(when a small periodic forcing term is present in addition to noise),
{\em coherence resonance} \cite{PiKu97}
(usually when the system is quiescent in the absence of noise),
and {\em self-induced stochastic resonance} \cite{HuDi93}
(involving relatively large noise that drives oscillations
of periods different from that of the deterministic system).

\begin{figure}[t!]
\begin{center}
\setlength{\unitlength}{1cm}
\begin{picture}(14.7,6.4)
\put(0,0){\includegraphics[width=7.2cm,height=6cm]{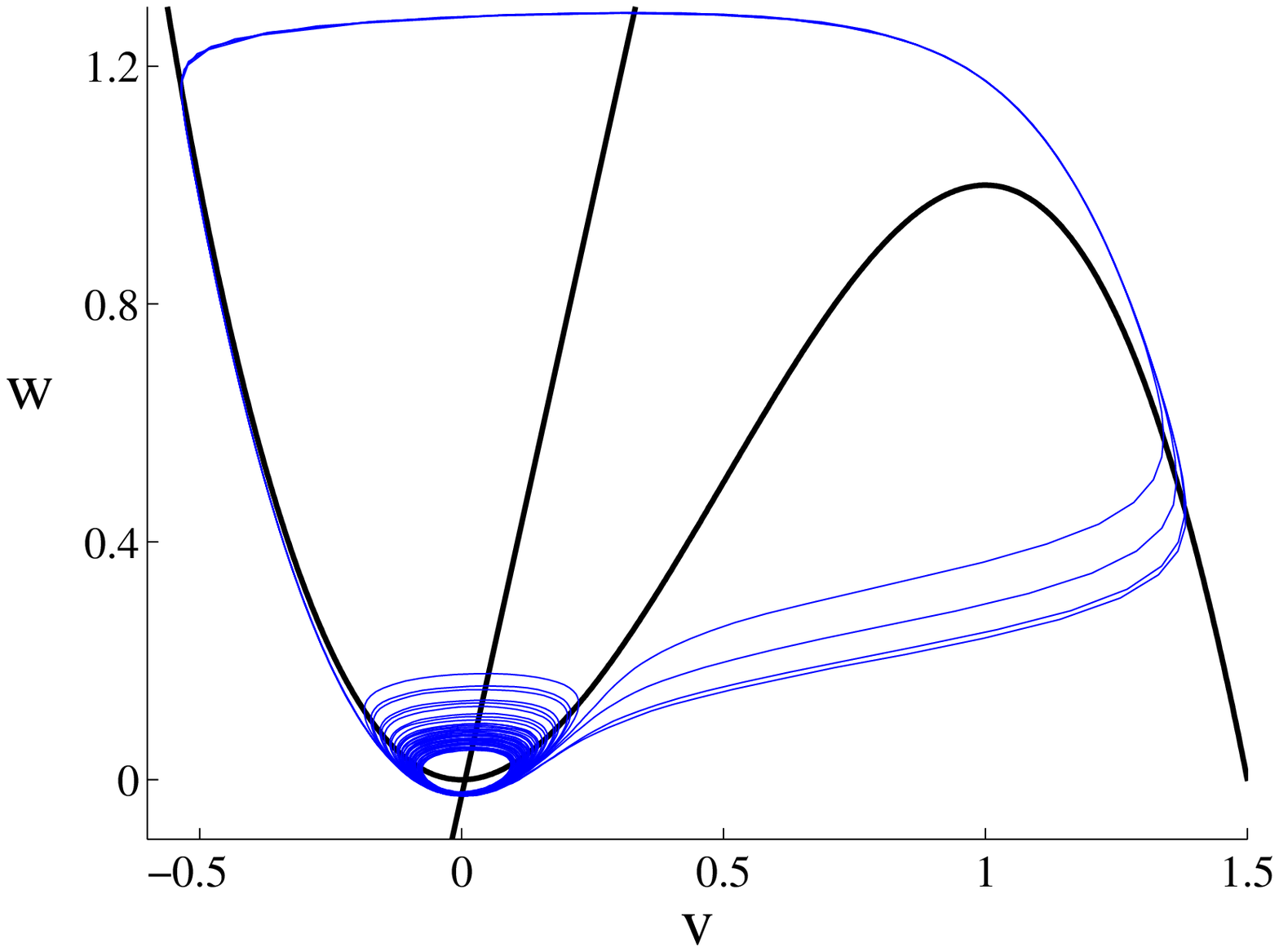}}
\put(7.7,0){\includegraphics[width=7.2cm,height=6cm]{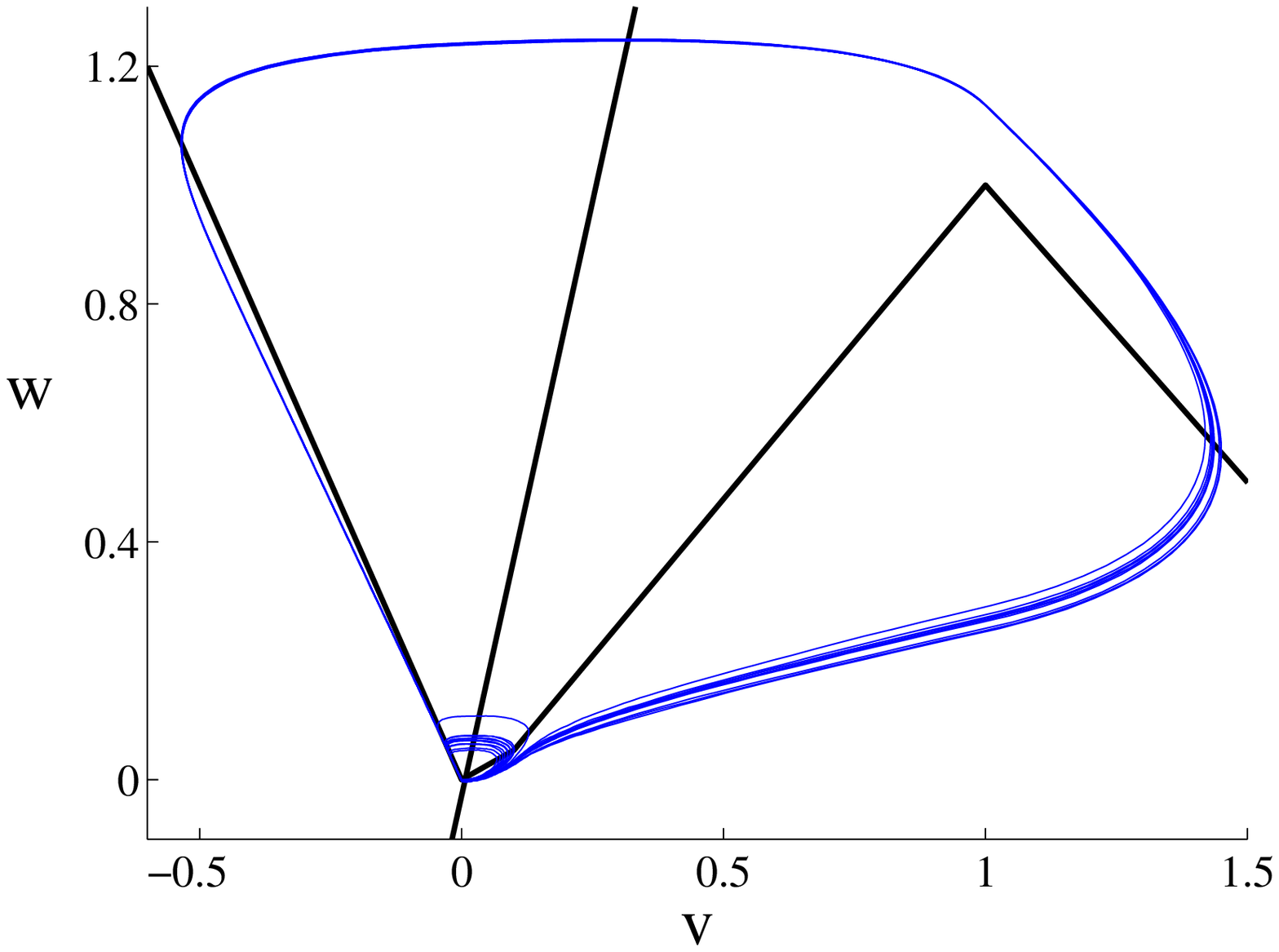}}
\put(1.2,6){\large \sf \bfseries A}
\put(8.9,6){\large \sf \bfseries B}
\put(3.4,6.3){smooth}
\put(11.1,6.3){PWL}
\end{picture}
\caption{
A trajectory of (\ref{eq:fhn}) with (\ref{eq:fcubic})
in panel A, and (\ref{eq:fhn}) with (\ref{eq:fpwl}) in panel B.
The parameter values used are the same as in
Fig.~\ref{fig:canardBifDiag}, with also $\lambda = 0.028$, $D = 0.0008$.
For these parameter values both the smooth and PWL models are
tuned to near the canard explosion.
\label{fig:ssEx}
}
\end{center}
\end{figure}

If $\lambda$ is tuned to values near the canard explosion,
in the presence of noise the system may exhibit both small amplitude and
large amplitude oscillations, i.e.~MMOs,
as shown in Figs.~\ref{fig:ssEx} and \ref{fig:timeSeries}.
Similar MMOs are described in \cite{DuMo08}
for very small noise by a careful choice of parameter values.
For a version of (\ref{eq:fhn})
that contains nonlinearity in the $w$ equation
to better mimic neural behaviour,
it has been observed that when
$\lambda$ is chosen to be just prior to the canard point
the frequency of relaxation oscillations increases
with noise amplitude \cite{MaNe01}.
Noise-induced MMOs have been described for three coupled
FHN systems near a canard \cite{LiWa07}.
A {\em signal-to-noise ratio} may be defined to quantitatively
determine dominant frequencies \cite{ZhHo05}.
Noise-induced MMOs may arise via
a different mechanism in the case that the Hopf bifurcation
is subcritical \cite{YuKu08}.
There are a variety of mechanisms for MMOs in three-dimensional systems
that we do not consider here, see for instance \cite{DeKr08,RoWe08}
and references in \cite{DeGu10}.

In this paper we study noise-driven MMOs in (\ref{eq:fhn}) with (\ref{eq:fpwl}).
We use analytical methods to identify parameter values
for which MMOs occur and describe the dependence of
each model parameter on MMOs.
For typical values of the noise amplitude, $D$,
MMOs occur over some interval of positive $\lambda$-values.
In order to find such intervals we determine exit distributions
for forward orbits of (\ref{eq:fhn}) with (\ref{eq:fpwl})
through various cross-sections of phase space.
The exit distributions allow us to deduce the amplitude of oscillations
and consequently find intervals of MMOs.
We show that MMOs are robust in the sense that
large variations in other model parameters
can have minimal effect on the width of the $\lambda$-intervals.

We note that the model we consider has additive noise in the $w$-equation only,
as in, for instance, \cite{MaNe01,PiKu97}.
This choice allows some simplifications in demonstrating the
analytical method, while still capturing qualitatively the behavior
that would be observed for more general additive noise.
Throughout the paper we indicate where this assumption allows some
simplification in the analysis, and we indicate the differences that
would need to be addressed for the case of noise also in the $v$-equation.

The remainder of the paper is organized as follows.
Section \ref{sec:DETER} briefly overviews PWL FHN models and
provides an analysis of (\ref{eq:fhn}) with (\ref{eq:fpwl})
in the absence of noise.
Here we explain the Hopf-like bifurcation at
$\lambda = 0$ that creates stable oscillations and
describe equations for $\lambda_{v_1}$ and $\lambda_1$,
Fig.~\ref{fig:lamEe}.
Calculations of exit distributions are detailed in \S\ref{sec:EXIT}.
Here we also describe the method by which we use these distributions
to find parameter values corresponding to MMOs.
Section \ref{sec:MMO} combines the analysis of the
previous sections to determine the effect of each model
parameter on MMOs.
Finally conclusions are presented in \S\ref{sec:CONC}.

\section{Properties of the deterministic system}
\label{sec:DETER}

\begin{figure}[ht]
\begin{center}
\includegraphics[width=9cm]{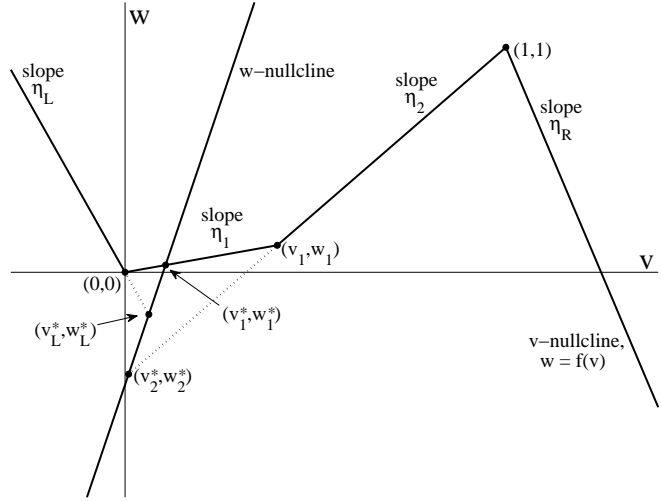}
\caption{
The nullclines of (\ref{eq:fhn}) with (\ref{eq:fpwl})
for small $\lambda > 0$.
Potential equilibria, $(v_j^*,w_j^*)$, lie at the intersection
of the nullclines.
If $\eta_j < \frac{\alpha}{\sigma}$ for every $j$,
then the system has a unique equilibrium
for all values of $\lambda$.
\label{fig:pp1}
}
\end{center}
\end{figure}

Analytical results may be derived for (\ref{eq:fhn})
when $f(v)$ is a PWL function.
Arguably the simplest continuous, PWL function that one can use for $f(v)$
consists of three line segments (one of them being the
straight connection between $(0,0)$ and $(1,1)$).
The FHN model with this function is well-studied \cite{LiSc00,Mc70},
refer to \cite{ItMu94} for the van der Pol system.
However, with this three-piece PWL function,
(\ref{eq:fhn}) does not exhibit a canard, as shown in \cite{ArOk97},
and so we do not consider it further.
Consequently, as in \cite{SeIn04,RoCo10}, we use two line segments
between $(0,0)$ and $(1,1)$ denoting the intermediate point by $(v_1,w_1)$
and the slopes by $\eta_j$, specifically (\ref{eq:fpwl}), as shown in Fig.~\ref{fig:pp1}.
If instead $f(v)$ contains multiple line segments left of $(0,0)$
such that the slopes of the two lines meeting at $(0,0)$ are $\pm \eta_1$,
multiple coexisting attractors commonly exist for small $\lambda$
which leads to complications that we do not study here.
For simplicity we do not consider $f(v)$
comprised of more than four line segments.
For canards in PWL FHN models with many segments we refer to reader 
to the recent work of Rotstein {\em et.~al.}~\cite{RoCo10}.

In the absence of noise (i.e.~when $D = 0$),
(\ref{eq:fhn}) with (\ref{eq:fpwl})
is a continuous, two-dimensional, PWL,
ordinary differential equation system:
\begin{equation}
\begin{split}
\dot{v} & = f(v) - w \;, \\
\dot{w} & = \varepsilon ( \alpha v - \sigma w - \lambda ) \;.
\end{split}
\label{eq:fhndet}
\end{equation}
The phase space, $\mathbb{R}^2$, is divided into four regions
\begin{equation}
\begin{split}
\mathcal{R}_L &= \{ (v,w) ~|~ v < 0, w \in \mathbb{R} \} \;, \\
\mathcal{R}_1 &= \{ (v,w) ~|~ 0 < v < v_1, w \in \mathbb{R} \} \;, \\
\mathcal{R}_2 &= \{ (v,w) ~|~ v_1 < v < 1, w \in \mathbb{R} \} \;, \\
\mathcal{R}_R &= \{ (v,w) ~|~ v > 1, w \in \mathbb{R} \} \;,
\end{split}
\label{eq:regions}
\end{equation}
by the three {\em switching manifolds}, $v=0$, $v=v_1$ and $v=1$,
on which the system is non-differentiable.

Each linear component of (\ref{eq:fhndet}) with (\ref{eq:fpwl})
has a unique equilibrium, $(v_j^*,w_j^*)$, see Fig.~\ref{fig:pp1}
(unless $\alpha = \sigma \eta_j$ in which case the relevant
$v$- and $w$-nullclines are parallel).
In the terminology of piecewise-smooth dynamical systems,
each $(v_j^*,w_j^*)$ is either {\em admissible}
(lies in the closure of $\mathcal{R}_j$)
or {\em virtual}
(lies outside the closure of $\mathcal{R}_j$).
The Jacobian, $A_j$, and the eigenvalues, $\rho_j$,
associated with each $(v_j^*,w_j^*)$ are 
\begin{eqnarray}
& A_j = \left[ \begin{array}{cc}
\eta_j & -1 \\
\varepsilon \alpha & -\varepsilon \sigma
\end{array} \right] \;, & \label{eq:Jac} \\
& \rho_j = \frac{1}{2} \left( \eta_j - \varepsilon \sigma \pm
\sqrt{ (\eta_j + \varepsilon \sigma)^2 - 4 \varepsilon \alpha} \right) \;. &
\label{eq:eigenvalues}
\end{eqnarray}
We assume
\begin{equation}
\eta_L < -\varepsilon \sigma - 2 \sqrt{\varepsilon \alpha} \;, \qquad
\varepsilon \sigma < \eta_1 < \frac{\alpha}{\sigma} \;,
\label{eq:etaAssump}
\end{equation}
such that $(v_L^*,w_L^*)$ is an attracting node
and $(v_1^*,w_1^*)$ is either a repelling node
or a repelling focus as determined by the sign of
$(\eta_1 + \varepsilon \sigma)^2 - 4 \varepsilon \alpha$.
The restriction (\ref{eq:etaAssump}) ensures that stable oscillations
are created at $\lambda = 0$, as shown below.

The bifurcation at $\lambda = 0$ that results from the interaction of
an equilibrium with the switching manifold, $v=0$,
is an example of
a {\em discontinuous bifurcation} \cite{DiBu08,LeNi04,Si10}.
Effectively, eigenvalues that determine the
stability of the admissible equilibrium change discontinuously
as the equilibrium crosses the switching manifold at $\lambda = 0$.
In general, a bifurcation is expected to occur
if one or more eigenvalues ``jump'' across the imaginary axis
at the crossing.
Such a bifurcation may be analogous to a smooth bifurcation
or it may be unique to piecewise-smooth systems \cite{LeNi04}.
For two-dimensional systems,
codimension-one, discontinuous bifurcations
involving a single smooth switching manifold have been
completely classified \cite{Si10,FrPo98}.

For the PWL system (\ref{eq:fhndet}) with (\ref{eq:fpwl}),
an attracting periodic orbit is born at the
discontinuous bifurcation, $\lambda = 0$.
The relative size of the periodic orbit for small $\lambda > 0$
is dependent upon whether the equilibrium, $(v_1^*,w_1^*)$,
is a node or a focus.
If $(v_L^*,w_L^*)$ is an attracting node
and $(v_1^*,w_1^*)$ is a repelling node,
invariant lines corresponding to eigenvectors
prevent the creation of a local periodic orbit
corresponding to a small oscillation \cite{Si10}.
The periodic orbit generated at $\lambda = 0$
has large amplitude (corresponding to a relaxation oscillation).
Specifically, as $\lambda \to 0^+$,
the maximum value of $v$ of the periodic orbit
limits on a value greater than $1$.

If instead $(v_1^*,w_1^*)$ is a repelling focus, then
the bifurcation is a discontinuous analogue of a Hopf bifurcation
in that a periodic orbit is created locally.
Unlike for a classical Hopf bifurcation,
the periodic orbit grows in size linearly
with respect to $\lambda$ (see Fig.~\ref{fig:canardBifDiag}-B)
which is typical for piecewise-smooth systems.

The value of $\varepsilon$ for which 
the square-root term in (\ref{eq:eigenvalues}) vanishes
is the critical value of $\varepsilon$
(see Fig.~\ref{fig:lamEe}-B) above which
the periodic orbit created at $\lambda = 0$ is small and below which
this orbit is large, and is given by
\begin{equation}
\varepsilon_{\rm crit} =
\frac{1}{\sigma^2} \left(
2 \alpha - \sigma \eta_1 - 2 \sqrt{\alpha(\alpha-\sigma \eta_1)} \right) \;.
\label{eq:eCrit}
\end{equation}
The curves $\lambda = \lambda_{v_1}(\varepsilon)$ and 
$\lambda = \lambda_1(\varepsilon)$, Fig.~\ref{fig:lamEe}-B,
which bound the region of medium oscillations,
emanate from $(\lambda,\varepsilon) = (0,\varepsilon_{\rm crit})$.
Since the underlying system is PWL,
we may obtain analytical expressions relating to these curves
by deriving the explicit solution to the flow of each linear component of
(\ref{eq:fhn}) with (\ref{eq:fpwl}).
We let $(v^{(j)}(t;v_0,w_0),w^{(j)}(t;v_0,w_0))$ denote the solution to the
linear component of (\ref{eq:fhndet}) with (\ref{eq:fpwl})
corresponding to $\mathcal{R}_j$, for an arbitrary initial condition, $(v_0,w_0)$.
For instance:
\begin{eqnarray}
\left[ \begin{array}{c}
v^{(1)}(t;v_0,w_0) \\
w^{(1)}(t;v_0,w_0) \end{array} \right] & = &
{\rm e}^{\frac{(\eta_1 - \varepsilon \sigma) t}{2}}
\left[ \begin{array}{c}
\cos(\omega_1 t) + \frac{\eta_1+\varepsilon \sigma}{2\omega_1} \sin(\omega_1 t) \\
\frac{\varepsilon \alpha}{\omega_1} \sin(\omega_1 t)
\end{array} \right. \nonumber \\
& & \left. \begin{array}{c}
-\frac{1}{\omega_1} \sin(\omega_1 t) \\
\cos(\omega_1 t) - \frac{\eta_1+\varepsilon \sigma}{2\omega_1} \sin(\omega_1 t)
\end{array} \right]
\left[ \begin{array}{c}
v_0 - v_1^* \\
w_0 - w_1^* \end{array} \right]
+ \left[ \begin{array}{c}
v_1^* \\
w_1^* \end{array} \right] \;,
\label{eq:explicit1}
\end{eqnarray}
which equals the solution to (\ref{eq:fhndet}) with (\ref{eq:fpwl})
for the same initial condition
whenever $(v^{(1)}(t),w^{(1)}(t))$ lies in the closure of $\mathcal{R}_1$
at all times between $0$ and $t$, and where
\begin{equation}
\omega_j = \frac{1}{2} \sqrt{ \left|
(\eta_j + \varepsilon \sigma)^2 - 4 \varepsilon \alpha \right| } \;.
\end{equation}

Unfortunately we cannot in general explicitly solve (\ref{eq:explicit1}) for $t$
(in particular solve: $v^{(1)}(t) = 0$ for $t$).
Consequently we are unable to
extract $\lambda_{v_1}$ or $\lambda_1$ explicitly in terms
of the parameters of the system.
For brevity we omit the details
and simply note that for the figures in this paper
we determine $\lambda_{v_1}$ and $\lambda_1$ by
numerically solving transcendental expressions.
This may be accomplished to any desired accuracy quickly
and does not require the use of a
differential equation solving method.

The following two approximations are used in the analysis of later sections.
For a wide range of parameter values the attracting periodic
orbit passes close to the origin.
If we approximate $\lambda_{v_1}$ by finding where
the next intersection of the forward orbit of $(0,0)$
with the $v$-nullcline is $(v_1,w_1)$,
then we obtain 
\begin{equation}
\lambda_{v_1} \approx \frac{\alpha v_1 - \sigma w_1}
{1 + {\rm e}^{\frac{(\eta_1 - \varepsilon \sigma) \pi}{2 \omega_1}}} \;,
\label{eq:lamApprox}
\end{equation}
which is particularly accurate for $\varepsilon \approx \varepsilon_{\rm crit}$,
as shown in Fig.~\ref{fig:lamEe}.

Second, the two eigenvalues associated
with $(v_L^*,w_L^*)$ (\ref{eq:eigenvalues})
are $\rho_{L,{\rm slow}} = O(\varepsilon)$ and
$\rho_{L,{\rm fast}} = \eta_L + O(\varepsilon)$.
Within $\mathcal{R}_L$, trajectories 
rapidly approach the associated slow eigenvector.
This eigenvector intersects the switching manifold, $v=0$, at
\begin{equation}
\hat{w}_L = \frac{\lambda \rho_{L,{\rm slow}}}{\alpha - \sigma \eta_L} \;.
\label{eq:wLhat}
\end{equation}
Consequently, trajectories such as large oscillations that spend a
relatively long period of continuous time in $\mathcal{R}_L$,
exit this region extremely close to the point $(0,\hat{w}_L)$.
This point is important below in the discussion of stochastic dynamics.
It is usually sufficient to approximate $\lambda_1$
by considering $(v^{(1)}(t;0,\hat{w}_L),w^{(1)}(t;0,\hat{w}_L))$
and the subsequent $(v^{(2)}(t),w^{(2)}(t))$
and finding the value of $\lambda$ where $(v^{(2)}(t),w^{(2)}(t))$ intersects $(1,1)$.
This is because $\lambda_1$ corresponds to the existence of a periodic orbit
with a maximum $v$-value of $1$,
which must intersect $(1,1)$,
see (\ref{eq:size}) and the surrounding discussion.

\section{Exit distributions}
\label{sec:EXIT}

To analyze noise-driven MMOs we consider solutions to (\ref{eq:fhn}) with (\ref{eq:fpwl})
in the presence of noise over long time frames
such that transient behaviour has decayed.
In this context we determine the fraction of oscillations that are small,
the fraction that are medium,
and the fraction that are large, referring to (\ref{eq:size}).
One method is to simply solve the system for a long time
and count the number of different oscillations.
This Monte-Carlo approach is useful for obtaining a basic
understanding of the system but poor for an
accurate quantitative analysis because the system must be
solved accurately for many parameter combinations
requiring considerable computation time.
Instead, since the system under consideration is PWL,
we are able to use exit distributions for the regions
(\ref{eq:regions}) to approximate these fractions.
This approach does not necessitate arbitrarily small $\varepsilon$.
In contrast, Muratov and Vanden-Eijnden \cite{MuVa08} applied
stochastic methods to (\ref{eq:fhn}) with (\ref{eq:fcubic})
by considering the system asymptotically
(i.e.~with arbitrarily small $\varepsilon$ and $\lambda$)
which essentially reduces the problem to one dimension.
In \cite{LiSc99}, the same system is considered but
in the limit $\varepsilon \to 0$ which also reduces
mathematical calculations to one dimension.

\begin{figure}[b!]
\begin{center}
\includegraphics[width=12cm]{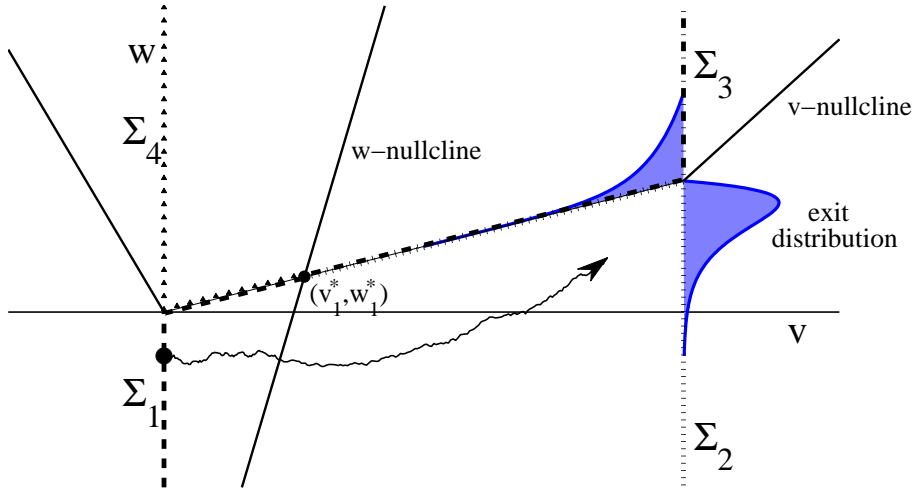}
\caption{
A sketch illustrating the exit distribution on $\Sigma_2$ (\ref{eq:Sigma})
for the forward evolution of a point on $\Sigma_1$.
For clarity each $\Sigma_j$ (\ref{eq:Sigma}) is drawn with
a different line type.
\label{fig:exitDist}
}
\end{center}
\end{figure}

Here we describe the exit distributions for forward orbits
of (\ref{eq:fhn}) with (\ref{eq:fpwl}) along various cross-sections
of phase space.
In the following section we use these exit distributions to
identify MMOs.
The four cross-sections we consider are:
\begin{equation}
\begin{split}
\Sigma_1 &= \{(0,w)~|~w < 0\} \cup \{(v,\eta_1 v)~|~0 \le v < v_1^*\} \;, \\
\Sigma_2 &= \{(v_1,w)~|~w < w_1\} \cup \{(v,\eta_1 v)~|~v_1^* < v \le v_1\} \;, \\
\Sigma_3 &= \{(v_1,w)~|~w > w_1\} \cup \{(v,\eta_1 v)~|~v_1^* < v \le v_1\} \;, \\
\Sigma_4 &= \{(0,w)~|~w > 0\} \cup \{(v,\eta_1 v)~|~0 \le v < v_1^*\} \;,
\end{split}
\label{eq:Sigma}
\end{equation}
as depicted in Fig.~\ref{fig:exitDist}.
We exclude the switching manifold, $v = 1$,
from calculations because large oscillations follow a
sufficiently predictable path back to $\mathcal{R}_L$ when $D \ll 1$.

One method for computing a first exit distribution is
to solve the Fokker-Planck equation
for the probability density of the process (\ref{eq:fhn}) with (\ref{eq:fpwl}) and
absorbing boundary conditions \cite{Sc80,Ga85}.
Integration of the solution to this boundary value problem
at the boundaries in an appropriate manner and over all positive time,
may yield the desired exit distribution.
However we dismiss this approach
as it necessitates extensive numerical computations,
in part because drift dominates the diffusion
which typically requires extra attention \cite{PaDe92,DoRu82,KnAn03}.
Instead we utilize the fact that within each region, $\mathcal{R}_j$,
the system is Ornstein-Uhlenbeck and,
ignoring switching manifolds, has a known explicit solution \cite{Ga85}.

The transitional probability density, $p_t^{(1)}$,
i.e.~${\rm Pr} \big( (v(t),w(t)) \in A ~\big|~
v(0) = v_0, w(0) = w_0, (v_0,w_0) \in \mathcal{R}_1 \big) =
\int\hspace{-1mm}\int\limits_{\hspace{-3mm}A}
p_t^{(1)}(v,w|v_0,w_0) \,dv\,dw$,
for the solution to
the Ornstein-Uhlenbeck process of $\mathcal{R}_1$,
i.e.~(\ref{eq:fhn}) with $f(v) = \eta_1 v$,
after a time $t$ is the Gaussian
\begin{equation}
p_t^{(1)}(v,w|v_0,w_0) = \frac{1}{2 \pi \sqrt{ \det(\Theta(t))}}
{\rm exp} \left( -\frac{1}{2}
\Delta z^{\sf T} \Theta(t)^{-1} \Delta z \right) \;,
\label{eq:gauss1}
\end{equation}
where
\begin{equation}
\Delta z(t;v_0,w_0) =
\left[ \begin{array}{c}
v - v^{(1)}(t;v_0,w_0) \\
w - w^{(1)}(t;v_0,w_0) \end{array} \right] \;.
\end{equation}
The mean, $(v^{(1)},w^{(1)})$, is the solution to the
system in absence of noise (\ref{eq:explicit1})
and $\Theta(t)$ is the covariance matrix given by:
\begin{equation}
\Theta(t) = D^2 \left[ \begin{array}{cc}
\int_0^t \left( {\rm e}^{A_1 s}_{12} \right)^2 \,ds &
\int_0^t {\rm e}^{A_1 s}_{12} {\rm e}^{A_1 s}_{22} \,ds \\
\int_0^t {\rm e}^{A_1 s}_{12} {\rm e}^{A_1 s}_{22} \,ds &
\int_0^t \left( {\rm e}^{A_1 s}_{22} \right)^2 \,ds \end{array} \right] 
= D^2 \left[ \begin{array}{cc}
\theta_{11}(t) & \theta_{12}(t) \\
\theta_{12}(t) & \theta_{22}(t)
\end{array} \right] \;,
\label{eq:Theta}
\end{equation}
where ${\rm e}^{A_1 s}_{ij}$ denotes the $(i,j)$-component of the
matrix exponential of $A_1 s$ (\ref{eq:Jac})
and we have introduced the $\theta_{ij}$ for convenience.
(Note that (\ref{eq:Theta}) would contain more terms
if (\ref{eq:fhn}) also included noise in the $v$ equation.)
The probability density (\ref{eq:gauss1})
obeys the Fokker-Planck equation
\begin{equation}
\frac{\partial p_t(v,w|v_0,w_0)}{\partial t} =
-\nabla \cdot J_t(v,w|v_0,w_0) \;,
\label{eq:FokkerPlanck}
\end{equation}
where
\begin{equation}
J_t(v,w|v_0,w_0) = \left[ \begin{array}{c}
(\eta_1 v - w) p_t \\
\varepsilon (\alpha v - \sigma w - \lambda) p_t -
\frac{D^2}{2}
\frac{\partial p_t}{\partial w} \end{array} \right] \;,
\label{eq:probCurrent}
\end{equation}
is the probability current \cite{Sc80,Ga85}.
By integrating (\ref{eq:FokkerPlanck}) and applying the divergence theorem,
it follows that the net flow of probability across, say,
the $v$-nullcline between $(v_1^*,w_1^*)$ and $(v_1,w_1)$,
is given by
\begin{equation}
\int_{v_1^*}^{v_1} n \cdot J_t(v,\eta_1 v|v_0,w_0) \,dv \;,
\nonumber
\end{equation}
where $n$ is the normal vector of the $v$-nullcline
pointing outwards \cite{Sc80,Ga85},
i.e.~here $n = \left[ -\frac{\eta_1}{1+\eta_1^2} ,
\frac{1}{1+\eta_1^2} \right]^{\sf T}$.
If trajectories were unable to cross the $v$-nullcline more than once,
then the integral
\begin{equation}
\int_0^{\infty} n \cdot J_t(v,\eta_1 v|v_0,w_0) \,dt \;,
\label{eq:exitDist}
\end{equation}
would be equal to the density of the first (and last)
exit points for escape from below the $v$-nullcline.
However, trajectories have multiple intersections with the $v$-nullcline
due to the presence of noise.
By considering two different time frames,
we now show that these multiple intersections have a negligible effect and
that (\ref{eq:exitDist}) represents an exit distribution suitable for our analysis.
Specifically we first show that the probability of return to the
$v$-nullcline after a short time is small.
Then we show that within this short time frame points of multiple
intersections are clustered.
Finally we show that for longer time intervals after an intersection with the $v$-nullcline,
trajectories are far from the nullcline, assuming small noise levels.

We first look at return times for the $v$-nullcline.
Closed form expressions for first passage problems of multi-dimensional
Ornstein-Uhlenbeck processes are not straight-forward
\cite{AlPa05,GrJa08}; for this reason we simplify to a one-dimensional problem.
Consider the forward orbit of a point
on the $v$-nullcline with $v_1^* < v < v_1$.
Using $y = w - \eta_1 v$ to represent the distance from the nullcline,
(\ref{eq:fhn}) with $f(v) = \eta_1 v$ may be written as
\begin{equation}
\begin{split}
dv & = -y~dt \;, \\
dy & = \big( (\eta_1 - \sigma \varepsilon) y +
(\alpha - \sigma \eta_1)(v - v_1^*) \varepsilon \big)~dt + D~dW \;,
\end{split}
\label{eq:fhnTransform1}
\end{equation}
where we have substituted 
$v_1^* = \frac{\lambda}{\alpha - \sigma \eta_1}$.
Intersections of the orbit with the nullcline are determined by
the $y$ equation of (\ref{eq:fhnTransform1})
which we conservatively reduce to
\begin{equation}
dy = c~dt + D~dW \;,
\label{eq:fhn1dApprox}
\end{equation}
where, for $y \ge 0$, the magnitude of the drift has a lower bound:
\begin{equation}
c \ge (\alpha - \sigma \eta_1)(v_{\rm min} - v_1^*) \varepsilon \;,
\nonumber
\end{equation}
assuming $v > v_{\rm min}$ for some $v_{\rm min} > v_1^*$.
For any $\delta > 0$, we are interested in
${\rm Pr} ( y(t)=0 {\rm ~for~some~} t \ge \delta ~|~ y(0)=0 )$,
i.e.~the probability that a solution to 
(\ref{eq:fhn1dApprox}) with $y(0) = 0$ satisfies $y(t) = 0$ at some $t \ge \delta$.
To calculate this probability we let $p(y,t)$ denote
the transitional probability density for (\ref{eq:fhn1dApprox}) with $y(0) = 0$,
and condition over the event that $y(\delta) = z$,
for all $z \in \mathbb{R}$:
\begin{equation}
{\rm Pr} \big( y(t)=0 {\rm ~for~some~} t \ge \delta ~\big|~ y(0)=0 \big) =
\int_{-\infty}^\infty p(z,\delta) \,
{\rm Pr} \big( y(t)=0 {\rm ~for~some~} t \ge \delta ~\big|~ y(\delta)=z \big) \,dz \;.
\nonumber
\end{equation}
Notice,
\begin{equation}
{\rm Pr} \big( y(t)=0 {\rm ~for~some~} t \ge \delta ~\big|~ y(\delta)=z \big)
= {\rm Pr} \big( y(t)=-z {\rm ~for~some~} t \ge 0 ~\big|~ y(0)=0 \big) \;,
\nonumber
\end{equation}
because (\ref{eq:fhn1dApprox}) has no explicit dependence on
$y$ and $t$.
This enables us to write
\begin{equation}
{\rm Pr} \big( y(t)=0 {\rm ~for~some~} t \ge \delta ~\big|~ y(0)=0 \big) =
\int_{-\infty}^\infty p(z,\delta) G(-z) \,dz 
\label{eq:appH2}
\end{equation}
where 
\begin{equation}
G(z) = {\rm Pr} \big( y(t)=z, {\rm ~for~some~} t \ge 0 ~\big|~ y(0)=0\big) \;.
\label{eq:appG}
\end{equation}
$G$ can be calculated from the density of the first hitting time
of $y(t)$ to $z$ (refer to \cite{Si51,Re01} for more details)
producing
\begin{equation}
G(z) = c \int_0^\infty p(z,t) \,dt \;.
\label{eq:appG2}
\end{equation}
By using (\ref{eq:appG2}) and
evaluating the integral on the right-hand side of (\ref{eq:appH2}), we obtain
\begin{equation}
{\rm Pr} \big( y(t)=0 {\rm ~for~some~} t \ge \delta ~\big|~ y(0)=0 \big)
= 1 - {\rm erf} \left(
\frac{c \sqrt{\delta}}{\sqrt{2} D}
\right) \;.
\end{equation}
For instance with $D = 0.0012$, $\varepsilon = 0.04$ and 
$(\alpha,\sigma) = (4,1)$,
whenever $v_{\rm min} - v_1^* > 0.025$ the probability of
return to the $v$-nullcline after a time of $\delta = 0.6$
is less than $1\%$. 

Second, for the system (\ref{eq:fhn}) with $f(v) = \eta_1 v$
we look at the distribution of future $v$-nullcline intersections
up to a time $\delta$.
The solution (\ref{eq:gauss1}) with $w_0 = \eta_1 v_0$ evaluated
on the $v$-nullcline and normalized is a
Gaussian with mean and variance:
\begin{eqnarray}
\tilde{v}(t) &=& \frac{-(\eta_1 \theta_{12} - \theta_{22}) v^{(1)}
+ (\eta_1 \theta_{11} - \theta_{12}) w^{(1)}}
{\theta_{22} - 2 \eta_1 \theta_{12} + \eta_1^2 \theta_{11}} \;, \label{eq:vTilde} \\
\tilde{\sigma}(t)^2 &=& \frac{\det{\Theta}}
{\theta_{22} - 2 \eta_1 \theta_{12} + \eta_1^2 \theta_{11}} \;, \label{eq:sigTilde}
\end{eqnarray}
respectively, where the $\theta_{ij}$ were defined in (\ref{eq:Theta}).
We observe that (\ref{eq:exitDist}) undergoes negligible
change when convolved by the Gaussian with 
(\ref{eq:vTilde}) and (\ref{eq:sigTilde}) evaluated at $t = \delta$.
For this reason we use (\ref{eq:exitDist})
to compute exit distributions on the $v$-nullclines.
The absence of noise in the $v$ equation of (\ref{eq:fhn})
ensures multiple rapid crossings through the switching manifolds are not permitted.
Consequently we use an integral similar to (\ref{eq:exitDist}) for exit distributions
across the other switching manifolds also.
This accounts for all components of each $\Sigma_j$ (\ref{eq:Sigma}).

\begin{figure}[t!]
\begin{center}
\setlength{\unitlength}{1cm}
\begin{picture}(13,6.3)
\put(0,0){\includegraphics[width=13cm,height=6.3cm]{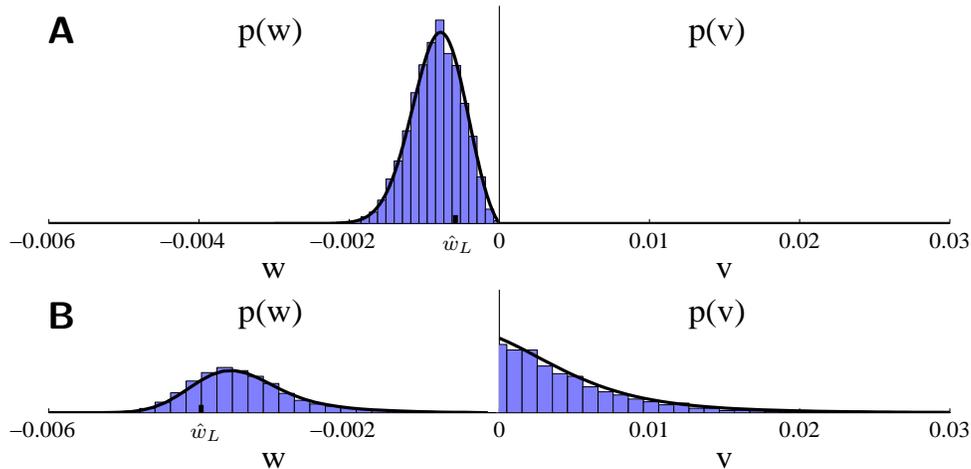}}
\put(.5,5.8){\large \sf \bfseries A}
\put(.5,2){\large \sf \bfseries B}
\put(2.4,.55){\scriptsize $\hat{w}_L$}
\put(5.75,3.05){\scriptsize $\hat{w}_L$}
\end{picture}
\caption{
Stationary densities of intersections of
(\ref{eq:fhn}) with (\ref{eq:fpwl}) on $\Sigma_1$
($v=0$ on the left, $w=\eta_1 v$ on the right).
The solid curves in panel B are computed using the iterative method
based on analytical expressions for the densities detailed in the text.
The curves in panel A could also be calculated by this iterative prodecure,
but instead it is more efficient to apply (\ref{eq:exitDist})
to the flow on the slow eigenvector of $\mathcal{R}_L$ with stationary variance.
This is because in panel A $\eta_1$ is relatively large
and oscillations enter $\mathcal{R}_L$ far from the origin
and so are strongly attracted to the slow eigenvector of $\mathcal{R}_L$.
The histograms are calculated from a single trajectory
of the system that was computed by numerical simulation
over a time period of $2 \times 10^5$.
The value of $\hat{w}_L$ (\ref{eq:wLhat}) is indicated in both panels.
In panel A, $w_1 = 0.05$ and $\lambda = 0.028$;
in panel B, $w_1 = 0.005$ and $\lambda = 0.19$.
The remaining parameter values are
$D = 0.008$,
$\varepsilon = 0.04$,
$v_1 = 0.1$,
$(\alpha,\sigma) = (4,1)$ and
$(\eta_L,\eta_R) = (-2,-1)$.
\label{fig:MCC}
}
\end{center}
\end{figure}

When the equilibrium, $(v_1^*,w_1^*)$, is admissible,
we expect the forward orbit of any point on $\Sigma_1$ to
escape the lower half of $\mathcal{R}_1$ (below the $v$-nullcline)
through $\Sigma_2$.
We calculate the exit distribution of the orbit through $\Sigma_2$
with (\ref{eq:exitDist}).
(Note, for simplicity we omit the
$\frac{\partial p}{\partial w}$ term in $J_t$ (\ref{eq:probCurrent})
when using (\ref{eq:exitDist}) because
it is dominated by the other terms in $J_t$.)
Using equally spaced data points and
performing this calculation repeatedly, we determine the
exit distribution on $\Sigma_2$
for any probability density of points on $\Sigma_1$.
From the exit distribution on $\Sigma_2$
we continue in a similar fashion and compute the
exit distributions on $\Sigma_3$, $\Sigma_4$ and lastly $\Sigma_1$.
Note these calculations use analytical expressions like (\ref{eq:exitDist})
and not Monte-Carlo simulations.
Numerically we observe that the iterative procedure of mapping
a distribution on $\Sigma_1$ to itself
(through $\Sigma_2$, $\Sigma_3$ and $\Sigma_4$) approaches the limiting
distribution of the intersection of an arbitrary forward orbit of the system
with $\Sigma_1$, Fig.~\ref{fig:MCC}.
We use the limiting distributions on $\Sigma_2$ and $\Sigma_3$
to calculate the probability that an arbitrary oscillation is
small, medium or large.
The results for a range of parameter values
are given in the next section.


\section{Mixed-mode oscillations}
\label{sec:MMO}

\begin{figure}[t!]
\begin{center}
\setlength{\unitlength}{1cm}
\begin{picture}(14.7,6.2)
\put(0,0){\includegraphics[width=7.2cm,height=6cm]{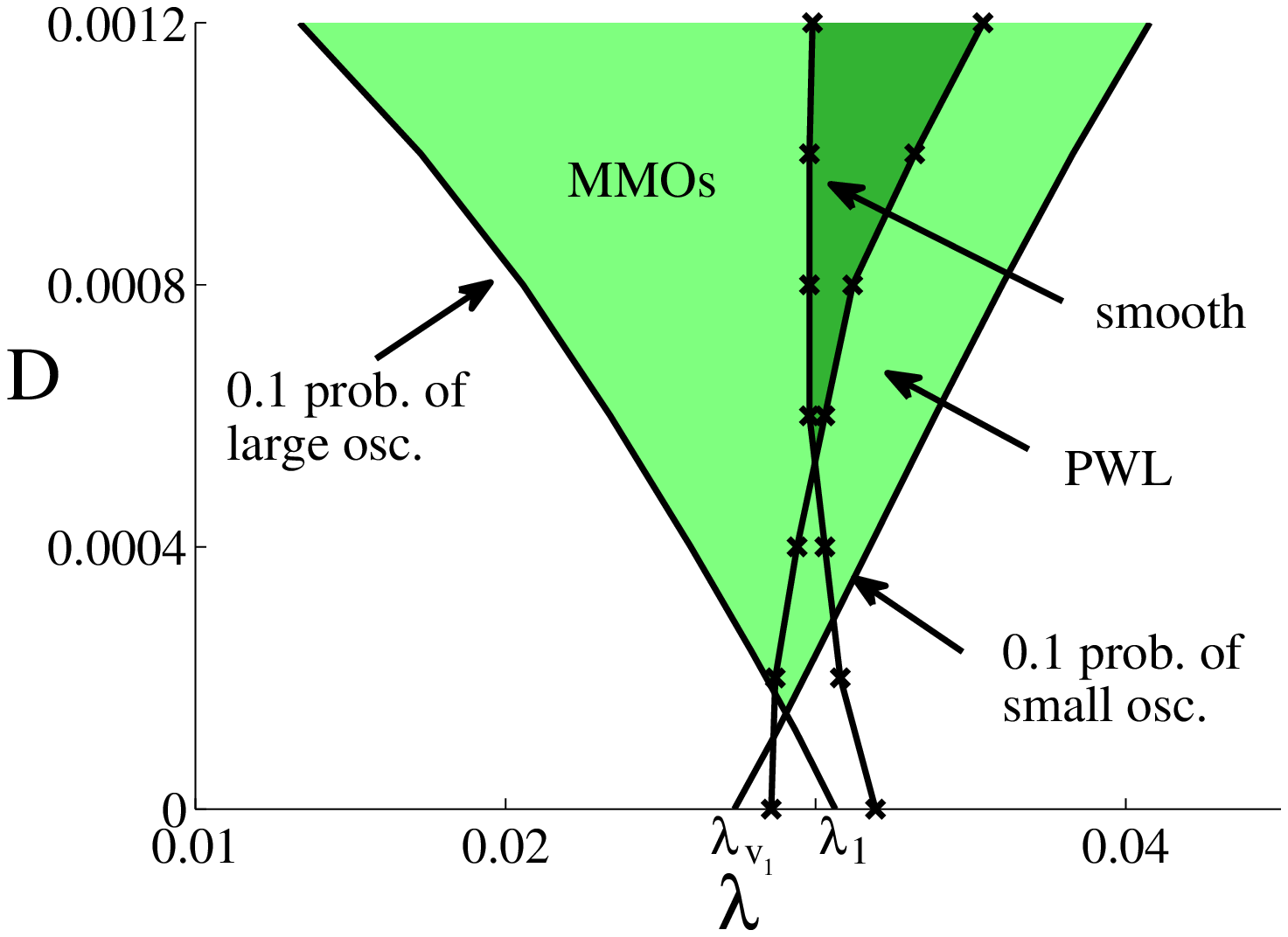}}
\put(7.7,0){\includegraphics[width=7.2cm,height=6cm]{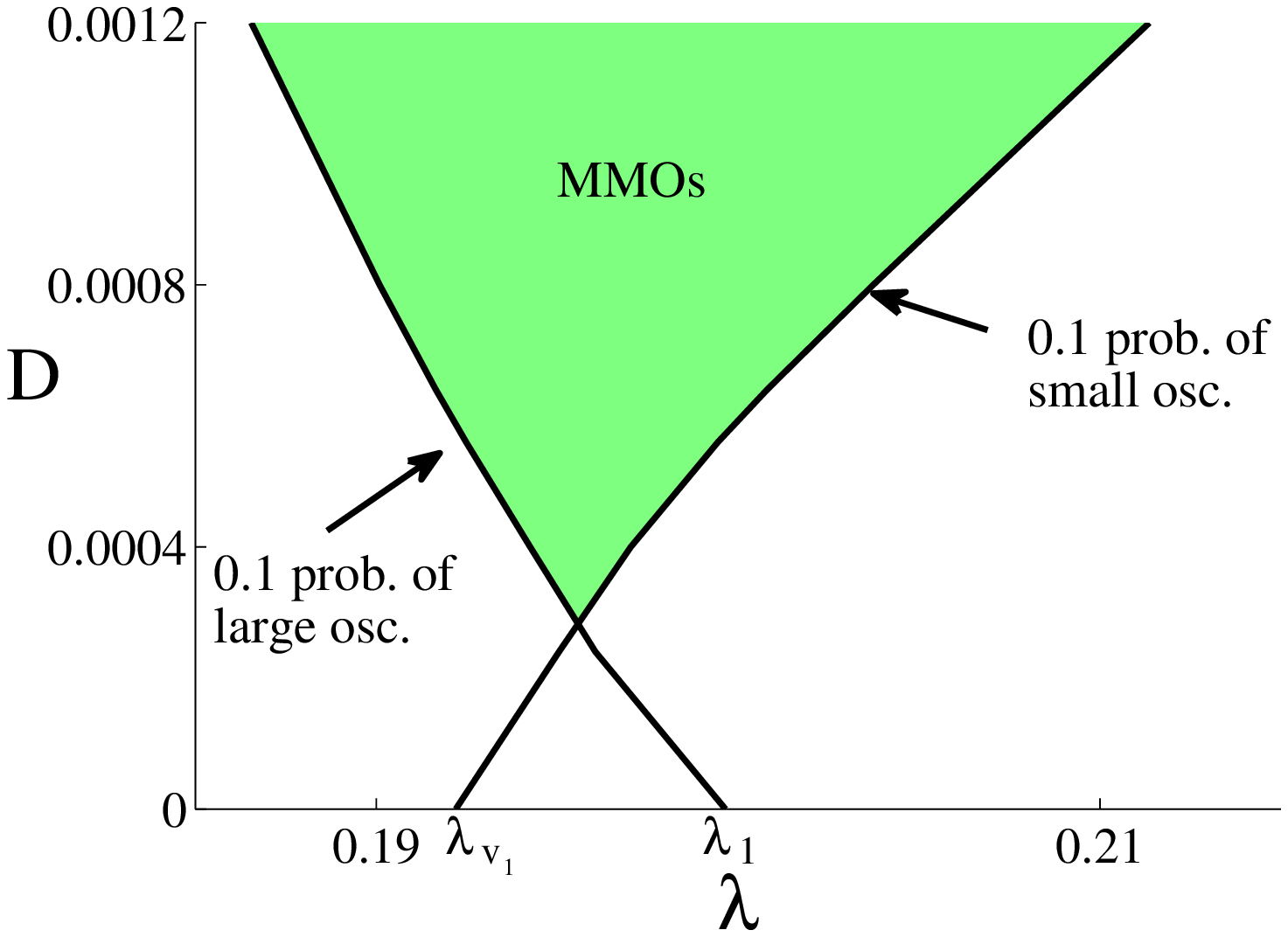}}
\put(0,5.2){\large \sf \bfseries A}
\put(7.7,5.2){\large \sf \bfseries B}
\end{picture}
\caption{
Regions of MMOs defined by where
at least $10\%$ of oscillations are small and
at least $10\%$ are large.
In panel A the parameter values used are
the same as in Fig.~\ref{fig:canardBifDiag}.
In panel B, $w_1 = 0.005$; the remaining parameter values are unchanged.
The region of MMOs for the smooth FHN model,
(\ref{eq:fhn}) with (\ref{eq:fcubic}), is superimposed in panel A.
\label{fig:lamXi}
}
\end{center}
\end{figure}

In order to understand MMOs quantitatively,
we say that (\ref{eq:fhn}) with (\ref{eq:fpwl})
exhibits MMOs whenever both small and large oscillations
occur at least $10\%$ of the time.
Specifically we find where exit densities corresponding to
small and large oscillations both integrate to a value greater than $0.1$.
Fig.~\ref{fig:lamXi} illustrates the dependence of MMOs
on the primary bifurcation parameter, $\lambda$,
and the noise amplitude, $D$.
Roughly the range of $\lambda$ values which permit MMOs increases with $D$.
This matches our intuition,
more noise allows for a wider variety of oscillations.
We compute Fig.~\ref{fig:lamXi} using the iterative scheme described in \S\ref{sec:EXIT};
Monte-Carlo simulations (not shown) give good agreement.

MMOs exist in a region bounded on the left by
the curve along which large oscillations occur $10\%$ of the time
and on the right by
the curve along which small oscillations occur $10\%$ of the time.
When $D = 0$ the former curve has the value $\lambda = \lambda_1$,
and the latter curve has the value $\lambda = \lambda_{v_1}$.
This is because the periodic orbit of
the system in the absence of noise, (\ref{eq:fhndet}),
changes from small to medium at $\lambda = \lambda_{v_1}$,
and from medium to large at $\lambda = \lambda_1$, \S\ref{sec:DETER}.

From Fig.~\ref{fig:lamXi},
we see that MMOs do not occur for arbitrarily small $D$
even near the canard explosion in contrast to what may be expected.
This is because for $D$ very small
and $\lambda_{v_1} < \lambda < \lambda_1$, medium oscillations dominate.
Medium oscillations occur less frequently with increasing $D$.
Note also that the MMO regions
appear relatively symmetric with respect to $\lambda$.

For the smooth system (\ref{eq:fhn}) with (\ref{eq:fcubic}),
we may roughly compute the region of MMOs, by the above definition,
from Monte-Carlo simulations, Fig.~\ref{fig:lamXi}-A.
We see that MMOs occur over a smaller parameter range for the
smooth version of the FHN model.
This distinction is possibly explained by Fig.~\ref{fig:ssEx}.
For the PWL model,
all oscillations (including small oscillations)
spend sufficient time in $\mathcal{R}_L$ to be drawn
into the slow eigenvector of this region.
Small and large oscillations are intertwined
in $\mathcal{R}_L$ on their approach to $\mathcal{R}_1$;
the amplitude of one oscillation is practically independent of
the previous oscillation.
(From a numerical viewpoint, in this situation fewer data points are required
than in general.)
In contrast, for the smooth system there is a significant distance
between small and large oscillations and therefore
more noise is required for, say, a large oscillation to
follow a small oscillation.

Noting this difference between the smooth and PWL models, we
considered another parameter range for the PWL system that has a
different exit distribution near the origin.
For small values of the slope, $\eta_1$, still respecting (\ref{eq:etaAssump}),
the MMOs may include small oscillations that do not enter $\mathcal{R}_L$,
so that the exit distribution across $\Sigma_1$ may be bimodal,
as shown in Fig.~\ref{fig:MCC}-B.
Here large oscillations intersect $\Sigma_1$ near $(0,\hat{w}_L)$ (\ref{eq:wLhat})
whereas the majority of small oscillations intersect $\Sigma_1$ on the $v$-nullcline.
We considered whether this type of bimodal exit distribution on
$\Sigma_1$ plays a role analogous to distance between
small and large oscillations in the
smooth model, but we did not see any evidence of this effect.
Specifically the MMO region, Fig.~\ref{fig:lamXi}-B,
has a similar size and shape to the region in Fig.~\ref{fig:lamXi}-A
for which the corresponding value of $\eta_1$ is an order of magnitude larger.

The boundaries of the MMO regions shown in Fig.~\ref{fig:lamXi} are relatively linear,
hence we perform an analytical calculation of the slopes at $D = 0$.
Let $s_{\rm small}$ [$s_{\rm large}$] denote the slope,
$\frac{d D}{d \lambda}$, at $D=0$, of the curve along which
$10\%$ of oscillations are small [large].

Let us begin with the curve along which exactly 10\% of oscillations are small.
This curve intersects $D = 0$ at $\lambda = \lambda_{v_1}$
at which the attracting periodic orbit created at $\lambda = 0$
intersects $w = \eta_1 v$ at $v = v_1$.
Here we can focus on small oscillations only,
so it suffices to consider the linear systems of $\mathcal{R}_L$ and $\mathcal{R}_1$,
i.e.~(\ref{eq:fhn}) with
\begin{equation}
f(v) = \left\{ \begin{array}{lc}
\eta_L v \;, & v \le 0 \\
\eta_1 v \;, & v > 0
\end{array} \right. \;.
\label{eq:L1only}
\end{equation}
As $\lambda$ is increased, the maximum $v$-value of
the deterministic periodic orbit increases at a rate, say, $\kappa_1$.
Due to linearity, this rate is given simply by
\begin{equation}
\kappa_1 = \frac{v_1}{\lambda_{v_1}} \;.
\label{eq:kappa1}
\end{equation}

\begin{figure}[b!]
\begin{center}
\setlength{\unitlength}{1cm}
\begin{picture}(9,5.4)
\put(0,0){\includegraphics[width=9cm,height=5.4cm]{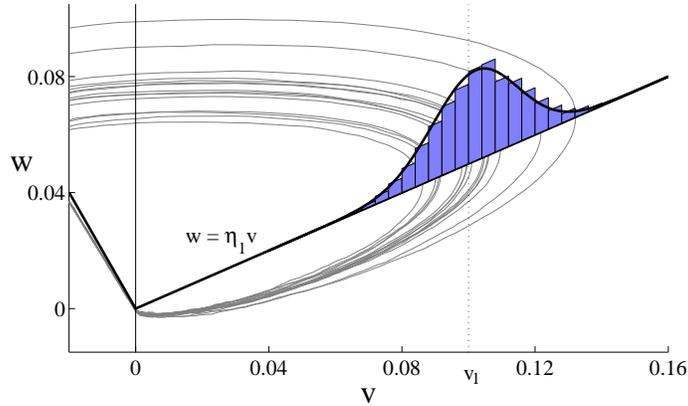}}
\end{picture}
\caption{
Intersections of (\ref{eq:fhn}) with (\ref{eq:L1only}) on $w = \eta_1 v$
using the same parameter values as Fig.~\ref{fig:canardBifDiag}
with also $\lambda = 0.028$, $D = 0.0004$.
The probability density curve is given by (\ref{eq:L1onlyDensity})
and the histogram is calculated from a single numerically computed trajectory
solved up to a time $2 \times 10^5$ (part of which is shown also).
\label{fig:vNullInt}
}
\end{center}
\end{figure}

When $\lambda = \lambda_{v_1}$,
the periodic orbit intersects $v=0$ at some point $(0,w_{\lambda_{v_1}})$
with $w_{\lambda_{v_1}} < 0$
and the line $w = \eta_1 v$ at $(v_1,w_1)$.
If we now consider small $D > 0$ but leave all other parameters unchanged,
over a long time frame trajectories intersect $v=0$
at points approximately normally distributed about $(0,w_{\lambda_{v_1}})$.
Since small oscillations neglect switching of (\ref{eq:fhn}) at $v=v_1$,
intersection points on $w = \eta_1 v$
are similarly approximately normally distributed about $(v_1,w_1)$,
as shown in Fig.~\ref{fig:vNullInt},
with a standard deviation of say, $\gamma_1 D$, where $\gamma_1$ is
a constant that we compute below.
The $v$-value of intersection points on $w = \eta_1 v$ then
have the distribution
$N(v_1 + \kappa_1(\lambda - \lambda_{v_1}),\gamma_1^2 D^2)$,
using (\ref{eq:kappa1}).
That is, if $q_{\rm small}$ denotes the probability density
for these $v$-values, then
\begin{equation}
q_{\rm small}(v) = \frac{1}{\sqrt{2 \pi} \gamma_1 D}
{\rm e}^{-\frac{1}{2 \gamma_1^2 D^2}
\left( v - v_1 - \kappa_1(\lambda-\lambda_{v_1}) \right)^2} \;.
\label{eq:L1onlyDensity}
\end{equation}
Then $10\%$ of oscillations are small when
\begin{equation}
\int_{-\infty}^{v_1} q_{\rm small}(v) \,dv =
\frac{1}{2} \left( 1 - {\rm erf} \left(
\frac{\kappa_1(\lambda-\lambda_{v_1})}{\sqrt{2} \gamma_1 D}
\right) \right) = 0.1 \;.
\end{equation}
By rearranging the previous equation we deduce that
the slope of the curve at $D = 0$ is 
\begin{equation}
s_{\rm small} =
\frac{d D}{d \lambda} =
\frac{\kappa_1}{\sqrt{2} \gamma_1 ~{\rm erf}^{-1}(0.8)} \;.
\label{eq:slopeSmall}
\end{equation}
From \S\ref{sec:DETER}, $\kappa_1$ may be accurately
calculated by solving transcendental equations.

We obtain a good approximation to $\gamma_1$ as follows.
Due to strong contraction in $\mathcal{R}_L$,
the distribution of points on $v=0$
has a standard deviation that is much
smaller than the standard deviation of points on $w = \eta_1 v$.
Than it is reasonable to approximate
the distribution on $v=0$ by the single value $w_{\lambda_{v_1}}$.
Then $q_{\rm small}$ is equivalent to the exit distribution along the
$v$-nullcline, thus by (\ref{eq:exitDist}),
\begin{equation}
q_{\rm small}(v) = \dot{w}(v,\eta_1 v)
\int_0^\infty p_t^{(1)}(v,\eta_1 v|0,w_{\lambda_{v_1}}) \,dt + O(D^2) \;,
\end{equation}
where $\dot{w}$ refers to (\ref{eq:fhndet}).
By (\ref{eq:gauss1}) and (\ref{eq:Theta}),
\begin{equation}
q_{\rm small}(v) = \frac{\dot{w}(v,\eta_1 v)}{2 \pi D^2}
\int_0^\infty \frac{1}{\sqrt{\theta_{11} \theta_{22} - \theta_{12}^2}} \;
{\rm e}^{-\frac{\phi(v,t)}{D^2}} \,dt + O(D^2) \;,
\label{eq:qSmall1}
\end{equation}
where
\begin{equation}
\phi(v,t) = \frac{1}{2 (\theta_{11} \theta_{22} - \theta_{12}^2)}
\left[ v-v^{(1)}, \eta_1 v - w^{(1)} \right]
\left[ \begin{array}{cc}
\theta_{22}(t) & -\theta_{12}(t) \\
-\theta_{12}(t) & \theta_{11}(t)
\end{array} \right]
\left[ \begin{array}{c}
v-v^{(1)} \\ \eta_1 v - w^{(1)}
\end{array} \right] \;.
\end{equation}
In the limit $D \to 0$,
the asymptotic approximation to integral in (\ref{eq:qSmall1}) is determined
from the main contribution of $\phi$, which is its maximum value;
formally this is achieved by Watson's lemma \cite{BeOr78}.
We omit the details of this calculation which produces
\begin{equation}
\gamma_1 = \sqrt{\theta_{1,1}} \;.
\end{equation}

We calculate the slope of the curve
on which $10\%$ of oscillations are large at $D = 0$
in a similar fashion.
Again we approximate the density of intersection points
on $v=0$ by a point mass but this time we use the value $\hat{w}_L$ (\ref{eq:wLhat})
and compute the density of intersection points
on the switching manifold, $v = v_1$.
For small $D$ this density is approximately Gaussian,
i.e.~$N(w_{\lambda_1} + \kappa_2(\lambda - \lambda_1),\gamma_2^2 D^2)$,
where when $\lambda = \lambda_1$
the deterministic trajectory passes through the points
$(0,\hat{w}_L)$ (or rather very near to this point),
$(v_1,w_{\lambda_1})$ and $(1,1)$.
If $q_{\rm large}$ denotes this probability density, then
\begin{equation}
q_{\rm large}(w) = \frac{1}{\sqrt{2 \pi} \gamma_2 D}
{\rm e}^{-\frac{1}{2 \gamma_2^2 D^2}
\left( w - w_{\lambda_1} - \kappa_2(\lambda-\lambda_1) \right)^2} \;.
\label{eq:L1onlyDensity2}
\end{equation}
The constant, $\kappa_2$,
may be computed from (\ref{eq:explicit1}) using the chain rule
for differentiation:
\begin{equation}
\kappa_2 = \frac{\partial w^{(1)}}{\partial \lambda}
= \frac{\partial w^{(1)}}{\partial t} \Bigg|_{t=t_{\rm int}}
\frac{\partial t_{\rm int}}{\partial \lambda} +
\frac{\partial w^{(1)}}{\partial w_0} \frac{\partial w_0}{\partial \lambda} +
\frac{\partial w^{(1)}}{\partial v_1^*} \frac{\partial v_1^*}{\partial \lambda} +
\frac{\partial w^{(1)}}{\partial w_1^*} \frac{\partial w_1^*}{\partial \lambda} \;.
\end{equation}
where $t_{\rm int}$ is the time taken for the trajectory
to go from $v = 0$ to $v = v_1$ and
\begin{equation}
\frac{\partial t_{\rm int}}{\partial \lambda} =
\left(
\frac{\partial v^{(1)}}{\partial w_0} \frac{\partial w_0}{\partial \lambda} +
\frac{\partial v^{(1)}}{\partial v_1^*} \frac{\partial v_1^*}{\partial \lambda} +
\frac{\partial v^{(1)}}{\partial w_1^*} \frac{\partial w_1^*}{\partial \lambda}
\right) \Bigg/
\frac{\partial v^{(1)}}{\partial t} \Bigg|_{t=t_{\rm int}} \;.
\end{equation}
By a calculation similar to that for $\gamma_1$ described above, we obtain
\begin{equation}
\gamma_2 =
\sqrt{\theta_{11} \left( \frac{\dot{w}}{\dot{v}} \right)^2 -
2 \theta_{12} \left( \frac{\dot{w}}{\dot{v}} \right) + \theta_{22}} \;,
\end{equation}
using (\ref{eq:fhndet}).

Unlike small oscillations, large oscillations traverse $\mathcal{R}_2$
and $\mathcal{R}_R$ so we must also consider the flow in these regions.
For $D = 0$ and $\lambda$ near $\lambda_1$,
we let $\hat{w}_2$ denote the first intersection of the
backwards orbit from $(1,1)$ with $v = v_1$
so that we may distinguish medium and large oscillations on $v=v_1$ when $D=0$.
We write 
\begin{equation}
\hat{w}_2 = w_{\lambda_1} + \kappa_3 ( \lambda - \lambda_1 ) \;,
\end{equation}
ignoring higher order terms and
where $\kappa_3$ may be calculated in a manner
similar to $\kappa_2$.
For small $D > 0$,
the forward orbit of any point $(v_1,w)$, with $w < w_1$,
has the probability, $p_{\rm large}(w;\lambda,D)$,
of undergoing a large, rather than medium, oscillation before returning to $\mathcal{R}_L$.
We find that for small $D$
there is a very sharp transition of $p_{\rm large}$ at $w = \hat{w}_2$ so
that it suffices to use the approximation
$p_{\rm large}(w) = H(\hat{w}_2 - w)$,
where $H(z) = 0$ for $z<0$ and $H(z) = 1$ for $z \ge 0$.
Consequently, $10\%$ of oscillations are large when
\begin{equation}
\int_{-\infty}^{w_{\lambda_1} + \kappa_3(\lambda-\lambda_1)} q_{\rm large}(w) \,dw =
\frac{1}{2} \left( 1 - {\rm erf} \left(
\frac{(\kappa_2-\kappa_3)(\lambda-\lambda_1)}{\sqrt{2} \gamma_2 D}
\right) \right) = 0.1 \;,
\end{equation}
where $q_{\rm large}$ is given by (\ref{eq:L1onlyDensity2}).
By rearranging this expression we arrive at
\begin{equation}
s_{\rm large} =
\frac{d D}{d \lambda} =
\frac{\kappa_2 - \kappa_3}{\sqrt{2} \gamma_2 ~{\rm erf}^{-1}(0.8)} \;.
\label{eq:slopeLarge}
\end{equation}
We have verified that computation of $s_{\rm small}$ and $s_{\rm large}$
by the expressions (\ref{eq:slopeSmall}) and (\ref{eq:slopeLarge})
matches with computation by exit distributions
(\ref{eq:exitDist}) as in Fig.~\ref{fig:lamXi}.

\begin{figure}[h!]
\begin{center}
\setlength{\unitlength}{1cm}
\begin{picture}(12,12.7)
\put(0,6.5){\includegraphics[width=12cm,height=6.2cm]{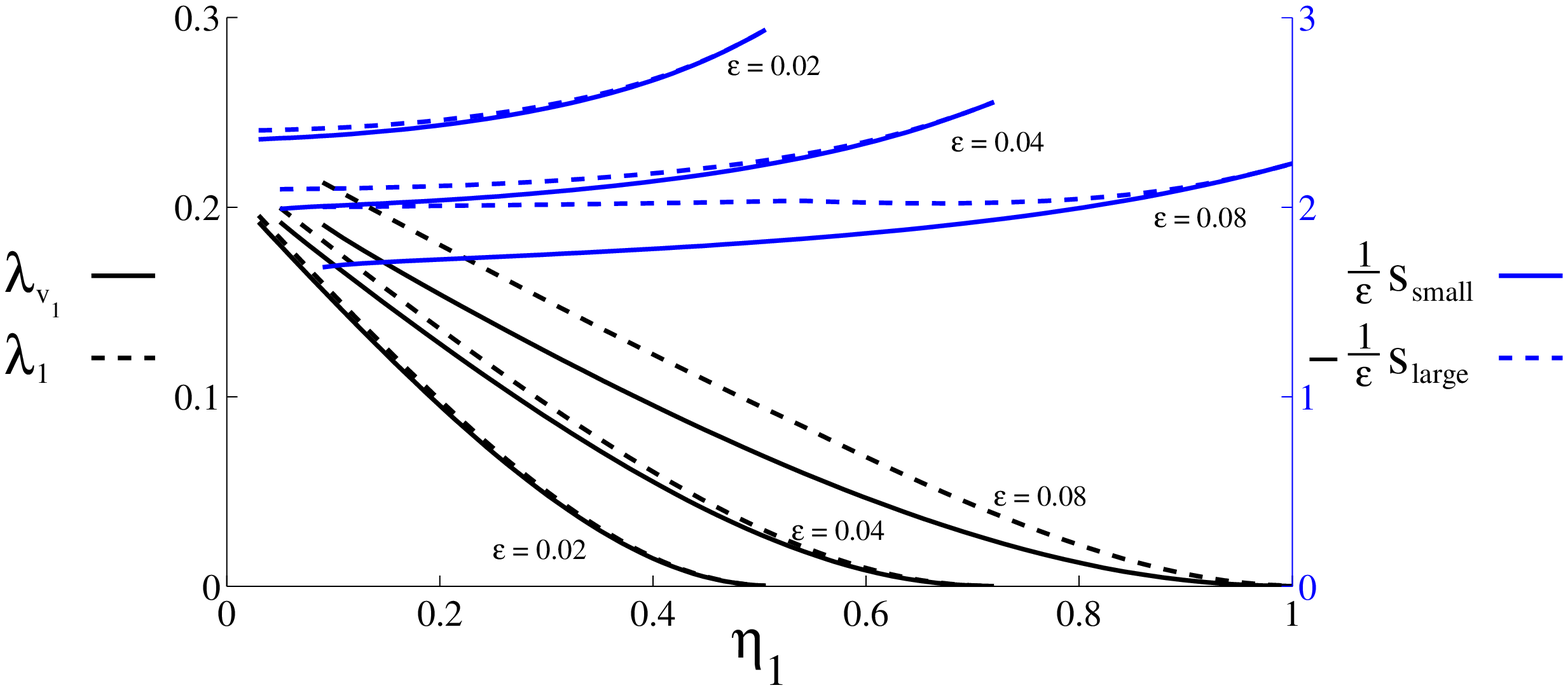}}
\put(0,0){\includegraphics[width=12cm,height=6.2cm]{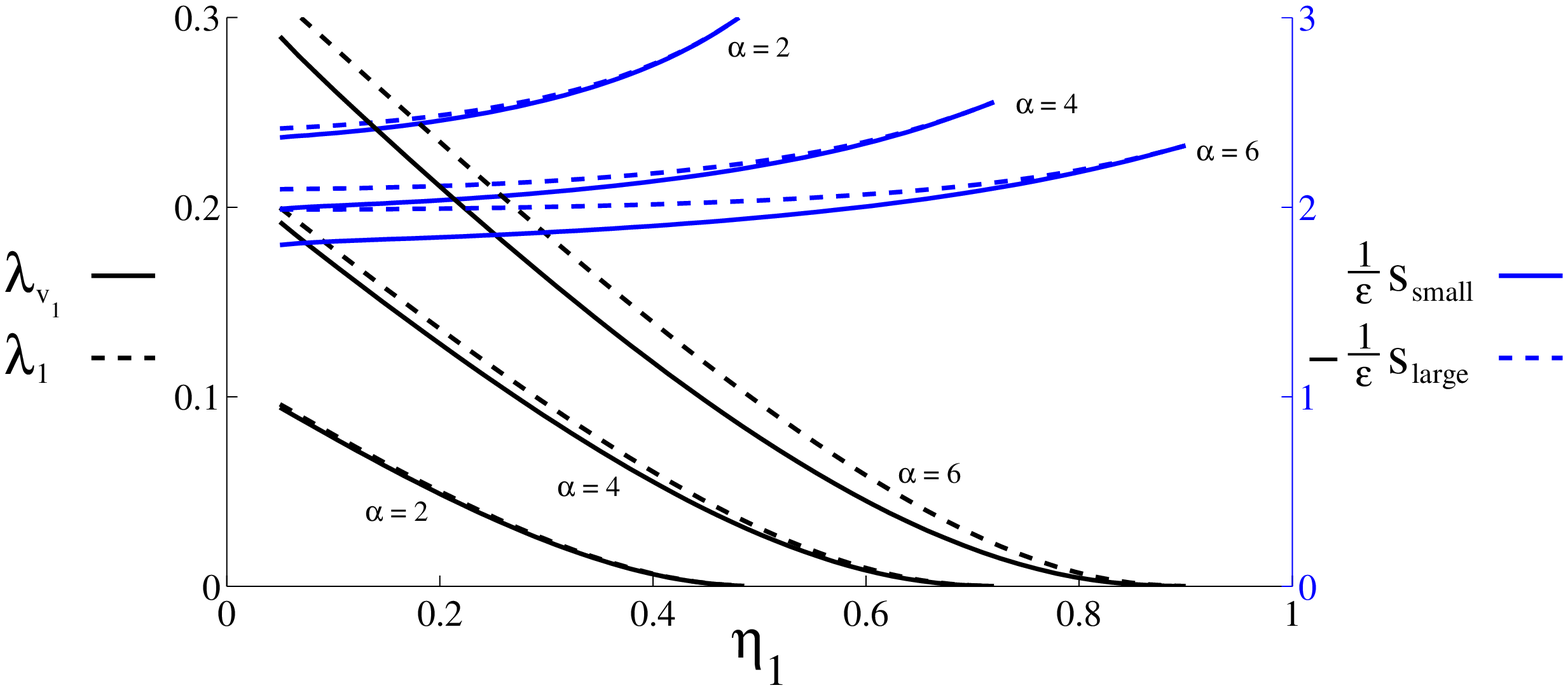}}
\put(.5,11.8){\large \sf \bfseries A}
\put(.5,5.3){\large \sf \bfseries B}
\end{picture}
\caption{
The dependence of $\lambda_{v_1}$, $\lambda_1$,
$s_{\rm small}$ and $s_{\rm large}$ on parameter values.
The lower curves are $\lambda_{v_1}$ and $\lambda_1$;
the upper curves are the slopes.
The $\lambda_{v_1}$ curves have end points at
$(\eta_1,\lambda) = (\varepsilon \sigma, \frac{1}{2}(\alpha - \varepsilon \sigma^2) v_1)$
and $(2\sqrt{\varepsilon \alpha} - \varepsilon \sigma,0)$.
In panel A, $\alpha = 4$.
In panel B, $\varepsilon = 0.04$.
In both panels, $\sigma = 1$, $v_1 = 0.1$, and
$(\eta_L,\eta_R) = (-2,-1)$.
We have scaled the slopes by $\frac{1}{\varepsilon}$ because the 
noise amplitude $D$ is not multiplied by $\varepsilon$ in (\ref{eq:fhn}).
\label{fig:slopeCalcs}
}
\end{center}
\end{figure}

The variation of the slopes (\ref{eq:slopeSmall}) and (\ref{eq:slopeLarge})
with respect to $\eta_1$, $\varepsilon$ and $\alpha$ is
shown in Fig.~\ref{fig:slopeCalcs}.
Both (\ref{eq:slopeSmall}) and (\ref{eq:slopeLarge})
approach zero as $\varepsilon \to 0$, but accounting for the fact
that the noise amplitude, $D$, is not multiplied by $\varepsilon$
in (\ref{eq:fhn}), the scaled values $\frac{1}{\varepsilon} s_{\rm small}$
and $\frac{1}{\varepsilon} s_{\rm large}$ vary relatively slightly.
With increasing $\eta_1$, $s_{\rm small}$ and $-s_{\rm large}$
increase slightly and approach the same value;
with increasing $\alpha$, $s_{\rm small}$ and $-s_{\rm large}$
decrease slightly.
Due to linearity, if $\eta_1$ is held constant and $v_1$ is increased,
$s_{\rm small}$ is unchanged.
We do not need to consider variation in $w_1$ and $\eta_2$
because these values may be written in terms of $\eta_1$ and $v_1$
(\ref{eq:eta1eta2}).
Additionally, the slopes are not strongly affected by $\eta_L$
(as long as $\eta_L \ll 0$) and $\eta_R$.
Therefore, for the most part, $\frac{1}{\varepsilon} s_{\rm small}$
and $-\frac{1}{\varepsilon} s_{\rm large}$
lie in, say, the interval $[1.8,2.5]$.
Hence for intermediate values of $\frac{D}{\varepsilon}$
although the interval of $\lambda$ values which permit MMOs
varies widely with the system parameters,
the width of this interval is robust with respect to parameter change.
For small values of $\frac{D}{\varepsilon}$, MMOs may not occur at all.
For very large values of $\frac{D}{\varepsilon}$,
multiple crossings on the $\Sigma_j$ may generate different results.

\section{Conclusions}
\label{sec:CONC}

We have studied MMOs in a PWL version of the FHN model,
(\ref{eq:fhn}) with (\ref{eq:fpwl}).
To obtain quantitative results we have defined oscillations
as small, medium, or large by the maximum $v$-value attained (\ref{eq:size}).
Furthermore we define MMOs by where at least $10\%$ of
oscillations are small and at least $10\%$ are large
(this approach may be applied to any preferred values of the percentages).
We incorporate noise additively in one equation for transparency of analysis.
Numerically we have observed that additive noise in both equations
yields similar mixed-mode dynamics.

The existence of a canard explosion in the system with no noise is a consequence
of incorporating a four-piece PWL function into the model.
We have introduced an analogy for canard points of smooth systems,
specifically we identify two values, $\lambda_{v_1}$ and $\lambda_1$,
at which the periodic orbit of the deterministic system
changes from small to medium, and from medium to large.
Both values increase with increasing $\alpha$, and $\varepsilon$,
and decrease with increasing $\eta_1$ as shown in Fig.~\ref{fig:slopeCalcs}.

Near the canard explosion noise drives MMOs.
The boundaries of the regions of MMOs shown in Fig.~\ref{fig:lamXi} are
approximately linear for small $D$ and we have calculated their slopes,
$\frac{d D}{d \lambda}$, at $D = 0$.
For small values of $D$ we have shown that the probability current
provides an approximation for exit distributions;
for large values of $D$ this approximation may no longer be valid.
Unless the noise amplitude is extremely small,
MMOs exist over some interval of $\lambda$ values.
We have illustrated that for constant
$\frac{D}{\varepsilon}$ typically the width of this interval
changes minimally with a relatively large variation in the values of the
other system parameters.

For the results in this paper we have used the value
$\eta_L = -2$ for the slope of the $v$-nullcline for $v < 0$.
With similar or more negative values of $\eta_L$ the system
exhibits the same qualitative behaviour.
However, with larger values of $\eta_L$, say $-1 \le \eta_L < 0$,
(\ref{eq:fhn}) with (\ref{eq:fpwl}) may have multiple attracting
solutions in the absence of noise.
The coexistence of attracting small and large periodic orbits,
naturally produces MMOs in the presence of noise
though this is via a different mechanism than the one studied here.

The FHN model of Makarov {\em et.~al.}~\cite{MaNe01},
which includes additional nonlinearity, readily exhibits MMOs.
It is possible that this addition strengthens the attraction
of the periodic orbit for $v < 0$, mimicking the slow eigenvector
discussed above and causing MMOs to be more robust.

\end{document}